\documentclass[12pt,a4paper]{article}
\usepackage{multirow}
\usepackage[dvips]{graphicx}
\usepackage{latexsym,amsmath,amsfonts,amscd}
\usepackage{epsfig}
\usepackage{indentfirst}
\usepackage{verbatim}
\usepackage{color}
\usepackage{colordvi}
\usepackage{tikz}
\usepackage{array}

\begin{document}

\baselineskip=2pc

\begin{center}
  {\Large \bf A  hybrid Hermite WENO scheme for hyperbolic conservation laws \footnote{The research is partly supported by Science Challenge Project, No. TZ2016002,  NSAF grant U1630247 and NSFC grant 11571290.} }
\end{center}

\centerline{Zhuang Zhao\footnote{School of Mathematical Sciences, Xiamen University,
Xiamen, Fujian 361005, P.R. China. E-mail: zzhao@stu.xmu.edu.cn.}, Yibing Chen\footnote{Institute of Applied Physics and Computational Mathematics, Beijing 100094, China. E-mail: chen\_yibing@iapcm.ac.cn.} and Jianxian
Qiu\footnote{School of Mathematical Sciences and Fujian Provincial
Key Laboratory of Mathematical Modeling and High-Performance
Scientific Computing, Xiamen University,
Xiamen, Fujian 361005, P.R. China. E-mail: jxqiu@xmu.edu.cn.}}

\vspace{.1in}

\baselineskip=1.8pc

\centerline{\bf Abstract}

\bigskip

\Red{In this paper, we propose a hybrid finite volume Hermite weighted essentially non-oscillatory (HWENO) scheme for solving one and two dimensional  hyperbolic conservation laws, which would be the fifth order accuracy in the one dimensional case, while is the  fourth order accuracy for two dimensional problems. The zeroth-order and the first-order moments are used in the spatial reconstruction, with total variation diminishing Runge-Kutta time discretization. Unlike the original HWENO schemes \cite{QSHw1,QSHw} using  different  stencils for spatial  discretization,  we borrow the thought of limiter for discontinuous Galerkin (DG) method to control the spurious oscillations, after this procedure, the scheme would avoid the oscillations by using HWENO reconstruction nearby discontinuities, and  using linear approximation straightforwardly in the smooth regions is to increase the efficiency of the scheme. Moreover, the scheme still keeps the compactness as only immediate neighbor information is needed in the reconstruction. A collection of benchmark numerical tests for one and two dimensional cases are performed to demonstrate the numerical accuracy, high resolution and robustness of the proposed scheme.}

\vfill {\bf Key Words:} Hermite WENO scheme, hyperbolic conservation
laws, discontinuous Galerkin method, limiter

{\bf AMS(MOS) subject classification:} 65M60, 35L65

\pagenumbering{arabic}

\newpage

\baselineskip=2pc

\section{Introduction}
\label{sec1}
\setcounter{equation}{0}
\setcounter{figure}{0}
\setcounter{table}{0}

In this paper, we design a hybrid  Hermite weighted essentially non-oscillatory (HWENO) scheme in the finite volume framework, which is the  fifth order accuracy in the one dimensional case and the fourth order accuracy for two dimensional problems. The HWENO scheme was derived from essentially non-oscillatory (ENO) and weighted essentially non-oscillatory (WENO) schemes, which have been widely applied for nonlinear hyperbolic conservation laws in recent decades. In 1985, Harten and Osher \cite{ho} constructed a weaker version of the total variation diminishing (TVD) criterion \cite{h}, which gave a framework for the reconstruction to design higher order ENO schemes. Then, Harten et al. \cite{heoc} developed the finite volume ENO schemes for solving one dimensional problems, in which they selected stencil adaptively in terms of the local smoothness, and in \cite{h1}, Harten extended the finite volume ENO schemes to two dimensional hyperbolic conservation laws. In 1994, the first WENO scheme was constructed by Liu, Osher and Chan \cite{loc}  mainly based on ENO scheme, where they used a nonlinear convex combination of all the candidate stencils to obtain higher order accuracy in smooth regions, and it was a third-order finite volume method in the one dimensional case. In 1996, Jiang and Shu \cite{js} proposed the third and fifth-order finite difference WENO schemes in multi-space dimensions, which gave a general framework to design the smoothness indicators and nonlinear weights. Ever since then, the WENO schemes have been further developed in the finite difference and finite volume frameworks presented in the literature \cite{hs,LPR,SHSN,ZQd,zqvo}, and more detailed review for WENO schemes can refer to \cite{s3}.

However, if we want to achieve higher order accuracy for WENO schemes, we need to enlarge the stencil for the spatial reconstruction. To make the stencil more compact, Qiu and Shu \cite{QSHw1,QSHw} developed the WENO methodology, which were first taken as limiters for Runge-Kutta discontinuous Galerkin methods, termed as Hermite WENO (HWENO) schemes. After this, many HWENO schemes were developed for solving hyperbolic conservations laws \cite{ZQHW,C0,HXA,TQ,LLQ,ZA,TLQ,CZQ,LLQ2}.
The HWENO schemes can achieve higher order accuracy than standard WENO schemes on the same stencils. On the other hand, to avoid spurious oscillations, all reconstructions of WENO and HWENO schemes  are based on local characteristic decompositions for the systems. Actually, the cost to compute the nonlinear weights and local characteristic decompositions is very high for WENO and HWENO schemes. To increase the efficiency, Pirozzoli \cite{SPir} designed an efficient hybrid compact-WENO scheme, on which they chose compact up-wind schemes in the smooth regions, while used WENO schemes in the discontinuous regions. Hill and Pullin \cite{Hdp} developed a hybrid scheme, which combined the tuned center-difference schemes with WENO schemes. The purpose was to expect the nonlinear weights would be achieved automatically in the smooth regions away from shocks, but a switching principle was still necessary. Next, Li and Qiu \cite{Glj} studied the hybrid WENO scheme using different switching principles \cite{TJS}, which illustrated that the troubled-cells indicator introduced by Krivodonova et al.\cite{LJJN} (KXRCF) has the ability to identify  the discontinuities well. Other different schemes introduced by \cite{ZQHD1,ZQHD2,ZZSQ,hby} for hyperbolic conservation laws also showed the good performances of the KXRCF troubled-cell indicator. In this paper, we would choose it as the indicator to identify the troubled-cell where the solutions may be discontinuous. The main idea for the hybrid WENO schemes \cite{Glj,ZZCQ} was that they both used non-linear WENO reconstruction near discontinuities, while employed upwind linear approximation directly in the smooth regions. We also notice that using upwind linear approximation in the smooth regions is one of  the choice, and the other  methods can also be used, such as B. Costa and W.S. Don used spectral method \cite{cd1} and central finite difference scheme \cite{cd2} in the smooth regions, which also can increase the efficiency obviously.

The hybrid HWENO scheme in this paper is different from the original HWENO schemes \cite{QSHw1,QSHw} as we take the thought of limiter for discontinuous Galerkin (DG) method to control the spurious oscillations. Since the solutions of nonlinear hyperbolic conservation laws often contain discontinuities, the derivative equations for the HWENO schemes need to deal with the derivatives or the first order moments, which would be relatively large nearby discontinuities. Therefore, the HWENO schemes listed by \cite{QSHw1,QSHw,ZQHW,TQ,LLQ,ZA,TLQ,CZQ,LLQ2} all used  different  stencils to  discretize the space for the original equations and the derivative equations, respectively. The variables of the derivative equations for the hybrid HWENO scheme are the first order moments, which also can be seen in the discontinuous Galerkin (DG) methods \cite{DG2,DG3,DG4,DG5} and other HWENO schemes \cite{ZQfeng,LWQdong,TLQ}. In one sense, these HWENO schemes can be seen as an extension by DG methods, and Dumbser et al. \cite{DBTM} gave a general and unified framework to define the numerical scheme extended by DG method, termed as $P_NP_M$ method, in which $P_N$ represents a piecewise polynomial of degree $N$ used as test functions in DG method and $P_M$ is a polynomial of degree $M$ reconstructed by the test functions of degree $N$  for  computing the numerical fluxes, and $M\geq N$. It is well known that DG methods  use limiters to modify the first order or higher order moments in the discontinuous regions, therefore, we adopt this thought  by adding a HWENO limiters \cite{QSHw1,QSHw,LBL} to modify the first order moments nearby discontinuities, and use HWENO procedure to reconstruct the point values on the interface of troubled cell.

The main procedures of the  hybrid HWENO scheme are given as follows. At first, we use the KXRCF troubled-cell indicator \cite{LJJN} to identify troubled cells, then, we  modify the first order moments in the troubled cells by the HWENO limiters \cite{QSHw1,QSHw,LBL}. After we modify the first order moments for all troubled cells in the computing domain, we would employ HWENO reconstruction at the points on the interface of the troubled cell, but use linear approximation at the internal points for spatial discretization, otherwise we directly use high order linear approximation.  For the systems, all HWENO reconstructions are based on local characteristic decompositions  to avoid spurious oscillations just like the classical WENO scheme \cite{js}. Compared with other HWENO schemes \cite{QSHw1,QSHw,ZQHW,C0,HXA,TQ,LLQ,ZA,TLQ,CZQ,LLQ2},  we borrow the idea of limiter for discontinuous Galerkin (DG) method to control the spurious oscillations, which would have two advantages. The one is to control the oscillations, and another one is to increase the efficiency for we directly use linear approximation in the smooth regions, where the limiter doesn't work under this circumstance. In short, the hybrid HWENO scheme avoids the spurious oscillations well and  has higher efficiency, while it still keeps the compactness as only immediate neighbor information is needed in the reconstruction.

The organization of the paper is as follows: in Section 2, we present the construction and implementation of  the finite volume hybrid HWENO scheme in the one and two  dimensional cases in detail. In Section 3, extensive numerical tests are performed to demonstrate the numerical accuracy, efficiency and robustness of the proposed scheme. Concluding remarks are given in Section 4.

\section{ Hybrid Hermite WENO scheme}
\label{sec2}
\setcounter{equation}{0}
\setcounter{figure}{0}
\setcounter{table}{0}
In this section, we will introduce the detailed implementation procedures of  the hybrid HWENO scheme for  one and two  dimensional hyperbolic conservational laws, which would be the fifth order accuracy in the one dimensional case, while is the  fourth order accuracy for two dimensional problems.

\subsection{One dimensional case}

We first consider one dimensional scalar hyperbolic conservation laws
  \begin{equation}
\label{EQ} \left\{
\begin{array}
{ll}
u_t+ f(u)_x=0, \\
u(x,0)=u_0(x). \\
\end{array}
\right.
\end{equation}
For simplicity, we divide the spatial domain with uniform meshes by $I_{i} =[x_{i-1/2},x_{i+1/2}]$, where the cell center is $x_i =\frac {x_{i-1/2}+x_{i+1/2}} 2$, and the mesh size is set as $ \Delta x =x_{i+1/2}- x_{i-1/2}$.

To design a HWENO scheme, we multiply the governing equation (\ref{EQ}) by $\frac{1}{\Delta x}$ and $\frac{x-x_i}{{(\Delta x)}^2}$, respectively, and  integrate  them over $I_i$, then, apply the integration by parts and use the numerical flux to approximate the values of the flux at the interface of $I_i$, finally, the semi-discrete finite volume HWENO scheme is
\begin{equation}
\label{odeH}
  \left \{
  \begin{aligned}
   \frac{d \overline u_i(t)}{dt} &=- \frac 1 {\Delta x} \left ( \hat f_{i+1/2}-\hat f_{i-1/2}\right ),\\
    \frac{d \overline v_i(t)}{dt}&=- \frac 1 {2\Delta x} \left ( \hat f_{i-1/2} +\hat f_{i+1/2}\right ) +\frac 1 {\Delta x} F_i(u),
  \end{aligned}
  \right.
\end{equation}
with initial conditions $\overline u_i(0) = \frac 1 {\Delta x} \int_{I_i} u_0(x) dx$ and $\overline v_i(0) = \frac 1 {\Delta x} \int_{I_i} u_0(x) \frac{x-x_i} { \Delta x} dx$, where $\overline u_i(t)$ is the cell average as $\frac 1 {\Delta x} \int_{I_i} u(x,t) dx$ and $\overline v_i(t)$ is the first order moment as $\frac 1 {\Delta x} \int_{I_i} u(x,t) \frac{x-x_i} { \Delta x}dx $. Here, $\hat f_{i+1/2}$ is the numerical flux which is the approximation to the values of  the flux $ f(u) $ at the interface point $x_{i+1/2}$ and $F_i(u)$ is the numerical integration for the flux $f(u)$. We adopt the  Lax-Friedrichs numerical flux method  to define the $\hat f_{i+1/2}$:
\begin{equation*}
\label{Nflux} \hat f_{i+1/2}= \frac 1 2 \left (f(u^-_{i+1/2})+f(u^+_{i+1/2}) \right) -\frac \alpha 2\left( u^+_{i+1/2}-u^-_{i+1/2} \right),
\end{equation*}
where $\alpha= \max_{u}|f'(u)|$. The numerical integration $F_i(u)$ is approximated by four-point Gauss-Lobatto quadrature formula, and the specific expression is given as follows,
\begin{equation*}
\label{Nfluxinte} F_i(u)=\frac 1 {\Delta x} \int_{I_i} f(u) dx \approx \sum_{l=1}^4 \omega_l f(u(x_l^G,t)).
\end{equation*}
Here, $ \omega_1=\omega_4=\frac 1 {12}$ and $ \omega_2=\omega_3=\frac 5 {12}$. The quadrature points on the cell $I_i$ are
\begin{equation*}
x_1^G=x_{i-1/2},\quad x_2^G=x_{i-\sqrt5/10}, \quad x_3^G=x_{i+\sqrt5/10}, \quad x_4^G=x_{i+1/2},
\end{equation*}
where $x_{i+a}$ is defined as $x_i +a \Delta x$.

\Red{The general frameworks for the hybrid HWENO scheme are given as follows. In Steps 1 and 2, we'll introduce the procedures of the spatial reconstruction for the semi-discrete scheme  (\ref{odeH}). In Step 3,  the equations (\ref{odeH}) is discretized in time by the third order TVD Runge-Kutta methodology \cite{so1}.}

\Red{\textbf{Step 1.} Identify the troubled-cell and  modify the first order moment in the troubled-cell.}

\Red{\textbf{Step 1.1.} Identify the troubled-cell.}

Troubled-cell means that the solution in the cell may be discontinuous, and in \cite{TJS}, Qiu and Shu investigated different troubled-cell indicators for Runge-Kutta discontinuous Galerkin methods. As suggested in \cite{TJS}, we use  KXRCF troubled-cell indicator  by Krivodonova et al. \cite{LJJN} to identify the discontinuities. We first divide the interface of the cell $I_i$ into  two parts $\partial I_i^-$ and $\partial I_i^+$, in which the flow is into ($\overrightarrow{v}\cdot\overrightarrow{n}<0$, $\overrightarrow{n}$ is the normal vector to $\partial I_i$) and  out ($\overrightarrow{v}\cdot\overrightarrow{n}>0$) of $I_i$, respectively. The cell $I_i$ is finally identified as a troubled cell, if:
\begin{equation}
\label{indicator}
 \frac{\left|\int_{\partial I_i^-}(u^h|_{I_i}-u^h|_{I_{n_i}})ds \right|} {h_i^{\frac{k+1}{2}}|\partial I_i^-|||u^h|_{I_i}||}>1,
\end{equation}
where $h_i$ is the radius of the circumscribed circle in the cell $I_i$, $I_{n_i}$ is the neighbor of $I_i$ on the side  of  $\partial I_i^-$, the norm is $L_\infty$ norm in the one dimensional case and $k$ is the degree of the polynomial $u_h$ approximating to $u(x)$, and we take $k=2$ in this paper. \Red{We should reconstruct the polynomial $u_h$ which is used only in the troubled-cell indicator (\ref{indicator}) to identify troubled cell,  not in the reconstruction procedure for solution.}

We use the information $\overline u_{i-1}$, $\overline u_{i}$, $\overline u_{i+1}$ and $\overline v_{i}$ to reconstruct a cubic polynomial $p_i^3(x)$ on the orthogonal basis function space $\left\{1, \frac {x-x_i}{\Delta x}, \left(\frac {x-x_i}{\Delta x}\right)^2-\frac 1{12}, \left(\frac {x-x_i}{\Delta x}\right)^3-\frac 3{20}\left(\frac {x-x_i}{\Delta x}\right) \right\}$, and the expressions is
\begin{equation*}
  p_i^3(x)=u_i^{(0)}+u_i^{(1)}\left(\frac {x-x_i}{\Delta x}\right)+u_i^{(2)}\left[ \left(\frac {x-x_i}{\Delta x}\right)^2-\frac 1{12} \right]+u_i^{(3)}\left[ \left(\frac {x-x_i}{\Delta x}\right)^3-\frac 3{20}\left(\frac {x-x_i}{\Delta x}\right) \right],
\end{equation*}
satisfying
\begin{equation*}
  \frac{1}{\Delta x} \int_{I_{i+j}} p_i^3(x)dx = \overline u_{i+j}, \quad j=-1,0,1, \quad \frac{1}{\Delta x} \int_{I_{i}} p_i^3(x)\frac{x-x_{i}} {\Delta x}dx = \overline v_{i}.
\end{equation*}
We have:
\begin{equation*}
u_i^{(0)}=\overline u_i,u_i^{(1)}=12\overline v_i,u_i^{(2)}=\frac12(\overline u_{i-1}-2\overline u_{i}+\overline u_{i+1}),u_i^{(3)}=-\frac 5{11}(\overline u_{i-1}+24\overline v_i -\overline u_{i+1}).
\end{equation*}
\Red{$u_h$ is taken as  $u_i^{(0)}+u_i^{(1)}\left(\frac {x-x_i}{\Delta x}\right)+u_i^{(2)}\left[ \left(\frac {x-x_i}{\Delta x}\right)^2-\frac 1{12} \right]$, which is the projection of $p_i^3(x)$ in the quadratic orthogonal function space $\left\{1,\frac {x-x_i}{\Delta x}, \left(\frac {x-x_i}{\Delta x}\right)^2-\frac 1{12} \right\}$.  Dropping the cubic term doesn't affect the accuracy for $u_h$  is used only in the troubled-cell indicator (\ref{indicator}) to identify troubled cell, we use next high order HWENO methodology to reconstruct the first order moment in the troubled cell.}

\Red{\textbf{Step 1.2.}  Modify the first order moment in the troubled-cell. }

If the cell $I_i$ is identified as a troubled cell, we would modify the first order moment $\overline v_i$.   The procedure to modify the first order moment is the same as  that HWENO limiter \cite{QSHw1}. At first, we give three  small stencils $S_1=\{I_{i-1},I_{i}\}$, $S_2=\{I_{i-1},I_{i},I_{i+1}\}$, $S_3=\{I_{i},I_{i+1}\}$,  and a large stencil $S_0=\{S_1,S_2,S_3\}=S_2$, then, we obtain three quadratic Hermite polynomials $p_1(x),p_2(x),p_3(x)$ on $S_1,S_2,S_3$, respectively, as
\begin{equation}
\begin{split}
\frac{1}{\Delta x} \int_{I_{i+j}} p_1(x)dx &= \overline u_{i+j}, \quad j=-1,0, \quad \frac{1}{\Delta x} \int_{I_{i-1}} p_1(x)\frac{x-x_{i-1}} {\Delta x}dx = \overline v_{i-1},\\
\frac{1}{\Delta x} \int_{I_{i+j}} p_2(x)dx &= \overline u_{i+j}, \quad j=-1,0,1, \\
\frac{1}{\Delta x} \int_{I_{i+j}} p_3(x)dx &= \overline u_{i+j}, \quad j=0,1, \quad \frac{1}{\Delta x} \int_{I_{i+1}} p_3(x)\frac{x-x_{i+1}} {\Delta x}dx = \overline v_{i+1},\\
\end{split}
\end{equation}
and get a quartic polynomial $p_0(x)$ on $S_0$, satisfying
\begin{equation}
 \frac{1}{\Delta x} \int_{I_{i+j}} p_0(x)dx = \overline u_{i+j}, j=-1,0,1, \quad \frac{1}{\Delta x} \int_{I_{i+j}} p_0(x)\frac{x-x_{i+j}} {\Delta x} dx = \overline v_{i+j}, j=-1,1.
\end{equation}
Then, we use these polynomials to reconstruct $\overline v_i$, and their explicit results based on the moments $\{\overline u_{i},\overline v_{i} \}_i$ are
\begin{equation*}
\begin{split}
\frac{1}{\Delta x} \int_{I_{i}} p_1(x)\frac{x-x_{i}} {\Delta x}dx &=\frac 1 6 \overline u_i-\frac 1 6 \overline u_{i-1}-\overline v_{i-1},\\
\frac{1}{\Delta x} \int_{I_{i}} p_2(x)\frac{x-x_{i}} {\Delta x}dx &= \frac 1 {24} \overline u_{i+1}-\frac 1 {24} \overline u_{i-1},\\
\frac{1}{\Delta x} \int_{I_{i}} p_3(x)\frac{x-x_{i}} {\Delta x}dx &= \frac 1 6 \overline u_{i+1}-\frac 1 6 \overline u_{i}-\overline v_{i+1},\\
\frac{1}{\Delta x} \int_{I_{i}} p_0(x)\frac{x-x_{i}} {\Delta x}dx &= \frac 5 {76} \overline u_{i+1}-\frac 5 {76} \overline u_{i-1}-\frac {11}{38}\overline v_{i-1}-\frac {11}{38}\overline v_{i+1}.\\
\end{split}
\end{equation*}
The linear weights $\gamma_1$, $\gamma_2$ and $\gamma_3$ can be obtained easily, just following as
\begin{equation*}
\frac{1}{\Delta x} \int_{I_{i}} p_0(x)\frac{x-x_{i}} {\Delta x}dx=\frac{1}{\Delta x} \sum_{n=1}^{3} \gamma_n \int_{I_{i}} p_n(x)\frac{x-x_{i}} {\Delta x}dx,
\end{equation*}
 which leads to $\gamma_1=\frac{11}{38}$, $\gamma_2=\frac{8}{19}$ and $\gamma_3=\frac{11}{38}$, then, we compute the smoothness indicators $\beta_n$, which measure how smooth the functions $p_n(x)$ in the target cell $I_i$, and we use the same definition for the smoothness indicators as in \cite{js,s3},
\begin{equation}
\label{GHYZ}
\beta_n=\sum_{\alpha=1}^r\int_{I_i}{\Delta x}^{2\alpha-1}(\frac{d ^\alpha
p_n(x)}{d x^\alpha})^2dx, \quad n=1,2,3.
\end{equation}
Here, $r=2$ is the degree of the polynomials $p_n(x)$, and their  explicit expressions are shown as
\begin{equation}
  \left \{
  \begin{aligned}
    \beta_1&=4(\overline u_{i-1}-\overline u_i+6 \overline v_{i-1})^2+\frac{13}{3}(\overline u_{i-1}-\overline u_i+12\overline v_{i-1})^2,\\
   \beta_2&=\frac{1}{4}( \overline u_{i-1}- \overline u_{i+1})^2+\frac{13}{12}(\overline u_{i-1}-2\overline u_i+\overline u_{i+1})^2, \\
    \beta_3&=4(\overline u_i-\overline u_{i+1}+6\overline v_{i+1})^2+\frac{13}{3}(\overline u_i-\overline u_{i+1}+12\overline v_{i+1})^2.
  \end{aligned}
  \right.
\end{equation}
Then the nonlinear weights are computed by:
\begin{equation*}
\omega_n=\frac{\bar\omega_n}{\sum_{l=1}^{3}\bar\omega_{l}},
\ \mbox{with} \ \bar\omega_{n}=\frac {\gamma_{n}}{(\beta_{n}+\varepsilon)^2}, \ n=1,2,3,
\end{equation*}
where $\varepsilon$ is a small positive number to avoid the denominator by zero, and we take $\varepsilon = 10^{-6}$ in this paper. Hence, the first order moment $\overline v_i$ is finally modified by
\begin{equation*}
  \overline v_i = \frac{1}{\Delta x} \sum_{n=1}^{3} \omega_n \int_{I_{i}} p_n(x)\frac{x-x_{i}} {\Delta x}dx.
\end{equation*}

\Red{\textbf{Step 2.} The reconstruction procedure for Gauss-Lobatto points values.}

\Red{In this subsection, we would give the details of reconstruction procedure for Gauss-Lobatto points values $u^\pm_{i\mp1/2}$ and $u_{i\pm\sqrt5/10}$ from $\{\overline u_{i},\overline v_{i} \}_i$.} Similarly to Step 1, we first give the stencils $S_1, S_2, S_3$ and $S_0$.  If one of the cells in stencil $S_0$ is identified as a troubled cell, we would use the HWENO method described in Step 2.1 to reconstruct the  $u^\pm_{i\mp1/2}$; otherwise we use the upwind linear approximation method described in Step 2.2 to reconstruct the  $u^\pm_{i\mp1/2}$.  And the  reconstruction of $u_{i\pm\sqrt5/10}$ is described in Step 2.3.

\Red{\textbf{Step 2.1.} The HWENO reconstruction for $u^-_{i+1/2}$.}

\Red{If one of the cells in stencil $S_0$ is identified as a troubled cell,  $u^-_{i+1/2}$ is reconstructed by  HWENO procedure.} We  reconstruct three cubic polynomials $p_1(x), p_2(x), p_3(x)$ on $S_1, S_2, S_3$, respectively, such that:
\begin{equation}
\label{HWEp1}
\begin{split}
\frac{1}{\Delta x} \int_{I_{i+j}} p_1(x)dx &= \overline u_{i+j}, \quad \frac{1}{\Delta x} \int_{I_{i+j}} p_1(x)\frac{x-x_{i+j}} {\Delta x}dx = \overline v_{i+j}, \quad j=-1,0,\\
\frac{1}{\Delta x} \int_{I_{i+j}} p_2(x)dx &= \overline u_{i+j}, \quad j=-1,0,1, \quad \frac{1}{\Delta x}  \int_{I_{i}} p_2(x)\frac{x-x_{i}} {\Delta x}dx = \overline v_{i},\\
\frac{1}{\Delta x} \int_{I_{i+j}} p_3(x)dx &= \overline u_{i+j}, \quad \frac{1}{\Delta x} \int_{I_{i+j}} p_3(x)\frac{x-x_{i+j}} {\Delta x}dx = \overline v_{i+j}, \quad j=0,1,\\
\end{split}
\end{equation}
and reconstruct a quintic polynomial $p_0(x)$ on $S_0$, as
\begin{equation}
\label{HWEp2}
 \frac{1}{\Delta x} \int_{I_{i+j}} p_0(x)dx = \overline u_{i+j}, \quad \frac{1}{\Delta x} \int_{I_{i+j}} p_0(x)\frac{x-x_{i+j}} {\Delta x} dx = \overline v_{i+j},\quad  j=-1,0,1.
\end{equation}
Based on  (\ref{HWEp1}) and (\ref{HWEp2}), we compute the approximations of $u^-_{i+1/2}$ by these polynomials at the point $x_{i+1/2}$, and their explicit expressions are
\begin{equation*}
\begin{split}
p_1(x_{i+1/2})&=\frac34\overline u_{i-1}+\frac14\overline u_i+\frac 7 2 \overline v_{i-1}+\frac{23}2 \overline v_i,\\
p_2(x_{i+1/2})&= \frac 2{33}\overline u_{i-1}+\frac 5 6 \overline u_i+ \frac 7 {66}  \overline u_{i+1}+\frac {60} {11} \overline v_i,\\
p_3(x_{i+1/2})&=\frac 1 2\overline u_{i}+\frac 1 2\overline u_{i+1}+2\overline v_{i}-2\overline v_{i+1},\\
p_0(x_{i+1/2})&=\frac {13} {108}\overline u_{i-1}+\frac 7{12}\overline u_i+\frac 8{27}\overline u_{i+1}+ \frac{25}{54}\overline v_{i-1}+\frac {241}{54}\overline v_i-\frac {28}{27}\overline v_{i+1}.
\end{split}
\end{equation*}
Then we get the linear weights $\gamma_1$, $\gamma_2$ and $\gamma_3$, according to
\begin{equation*}
 p_0(x_{i+1/2})=\sum_{n=1}^{3} \gamma_n p_n(x_{i+1/2}),
\end{equation*}
and we have  $\gamma_1=\frac{25}{189}$, $\gamma_2=\frac{22}{63}$ and $\gamma_3=\frac{14}{27}$. Then, we compute the smoothness indicators $\beta_n$, which measure how smooth the functions $p_n(x)$ in the cell $I_i$.  Again, we use the formula (\ref{GHYZ}) to compute the smoothness indicators,  where $r=3$, and we have the  explicit expressions for the smoothness indicators:
\begin{equation}
\label{GHYZP}
\left \{
\begin{aligned}
\beta_1&=\frac 1 {16} (\overline u_{i-1}-\overline u_i+6\overline v_{i-1}+54\overline v_{i} )^2+\frac {13}{48}(15\overline u_{i-1}-15\overline u_i+66\overline v_{i-1}+114\overline v_{i} )^2+\\
       &\quad \frac {3905} {16} (\overline u_{i-1}-\overline u_i+6 \overline v_{i-1}+6\overline v_{i} )^2,\\
\beta_2&=\frac 1 {484}(\overline u_{i-1}-\overline u_{i+1}-240\overline v_{i} )^2+\frac {13}{12}(\overline u_{i-1}-2\overline u_i+\overline u_{i+1} )^2+\\
       &\quad \ \frac {355} {44} (\overline u_{i-1}-\overline u_{i+1}+24\overline v_{i})^2, \\
\beta_3&=\frac 1 {16} (\overline u_{i}-\overline u_{i+1}+54\overline v_{i}+6\overline v_{i+1} )^2+\frac {13}{48}(15\overline u_{i}-15\overline u_{i+1}+114\overline v_{i}+66\overline v_{i+1} )^2+\\
       &\quad  \frac {3905} {16} (\overline u_{i}-\overline u_{i+1}+6 \overline v_{i}+6\overline v_{i+1} )^2,\\
\end{aligned}
\right.
\end{equation}
and the nonlinear weights are computed by:
\begin{equation*}
\omega_n=\frac{\bar\omega_n}{\sum_{l=1}^{3}\bar\omega_{l}},
\ \mbox{with} \ \bar\omega_{n}=\frac {\gamma_{n}}{(\beta_{n}+\varepsilon)^2}, \ n=1,2,3.
\end{equation*}
Here $\varepsilon$ is a small positive number taken as $10^{-6}$. Hence, the final value of $u^-_{i+1/2}$ is reconstructed by
\begin{equation*}
  u^-_{i+1/2} = \sum_{n=1}^{3} \omega_n p_n(x_{i+1/2}).
\end{equation*}

The reconstruction to $u^+_{i-1/2}$ is mirror symmetric with respect to $x_i$ of the above procedure.

\Red{\textbf{Step 2.2.} The linear approximation for $ u^\mp_{i\pm1/2}$.}

\Red{If nether cell in $S_0$ is identified as troubled cell, then we will use upwind linear reconstruction for $ u^\mp_{i\pm1/2}$, that is we use  $p_0(x)$ to approximate $u$ directly, and we have:
\begin{equation*}
u^+_{i-1/2} = p_0(x_{i-1/2})=\frac 8{27}\overline u_{i-1}+\frac 7{12}\overline u_i+\frac {13} {108}\overline u_{i+1}+\frac {28}{27} \overline v_{i-1}-\frac {241}{54}\overline v_i-\frac{25}{54}\overline v_{i+1},
\end{equation*}
and
\begin{equation*}
u^-_{i+1/2} =  p_0(x_{i+1/2})=\frac {13} {108}\overline u_{i-1}+\frac 7{12}\overline u_i+\frac 8{27}\overline u_{i+1}+ \frac{25}{54}\overline v_{i-1}+\frac {241}{54}\overline v_i-\frac {28}{27}\overline v_{i+1}.
\end{equation*}}

\Red{\textbf{Step 2.3.} The linear reconstruction for $u_{i\pm\sqrt5/10}$.}

\Red{We would like to use linear reconstruction for $u_{i\pm\sqrt5/10}$ in all cells, then, $u_{i\pm\sqrt5/10}$ are finally approximated by the following expressions, respectively,
\begin{equation*}
\begin{split}
u_{i-\sqrt5/10} = p_0(x_{i-\sqrt5/10})&=-(\frac{101}{5400}\sqrt5+\frac{1}{24})\overline u_{i-1}+\frac{13}{12}\overline u_i+(\frac{101}{5400}\sqrt5-\frac{1}{24})\overline u_{i+1}-\\
&\quad (\frac{3}{20}+\frac{841}{13500}\sqrt5)\overline v_{i-1}-\frac{10289}{6750}\sqrt5\overline v_i+(\frac{3}{20}-\frac{841}{13500}\sqrt5)\overline v_{i+1},
\end{split}
\end{equation*}
and
\begin{equation*}
\begin{split}
u_{i+\sqrt5/10} = p_0(x_{i+\sqrt5/10})&=(\frac{101}{5400}\sqrt5-\frac{1}{24})\overline u_{i-1}+\frac{13}{12}\overline u_i-(\frac{101}{5400}\sqrt5+\frac{1}{24})\overline u_{i+1}+\\
&\quad (\frac{841}{13500}\sqrt5-\frac{3}{20})\overline v_{i-1}+\frac{10289}{6750}\sqrt5\overline v_i+(\frac{3}{20}+\frac{841}{13500}\sqrt5)\overline v_{i+1}.
\end{split}
\end{equation*}}

\textbf{Step 3.} When we have finished the spatial discretization following Steps 1 and 2, the semi-discrete schemes (\ref{odeH}) are discretized in time by the third order TVD Runge-Kutta method \cite{so1}:
\begin{eqnarray}
\label{RK}\left \{
\begin{array}{lll}
     u^{(1)} & = & u^n + \Delta t L(u^n),\\
     u^{(2)} & = & \frac 3 4u^n + \frac 1 4u^{(1)}+\frac 1 4\Delta t L(u^{(1)}),\\
     u^{(n+1)} & = &\frac 1 3 u^n +  \frac 2 3u^{(2)} +\frac 2 3\Delta t L(u^{(2)}).
\end{array}
\right.
\end{eqnarray}

{\bf \em Remark 1:} For one dimensional scalar equation, the solution $u$ is taken as the indicator variable and $\overrightarrow v$ is defined as $f'(u)$; while for one dimensional Euler equations, the density $\rho$ and the energy $E$ are taken as the indicator variables, respectively, and $\overrightarrow v$ is set as the velocity $\mu$ of the fluid.

{\bf \em Remark 2:}  For the systems, such as the one dimensional compressible Euler equations,   in order to achieve better qualities at the price of more complicated computations, the HWENO approximation is always used with a local characteristic field decomposition seen in e.g. \cite{s2,s3} for details, while the  linear approximation is used in component by component.

\subsection{Two dimensional case}

Similarly to one dimensional case, we first consider two dimensional scalar hyperbolic conservation laws
  \begin{equation}
\label{EQ2} \left\{
\begin{array}
{ll}
u_t+ f(u)_x+g(u)_y=0, \\
u(x,y,0)=u_0(x,y). \\
\end{array}
\right.
\end{equation}
For simplicity of presentation, the computing domain is divided by uniform meshes. The mesh sizes are $ \Delta x =x_{i+1/2}- x_{i-1/2}$ in the $x$ direction and  $ \Delta y =y_{j+1/2}- y_{j-1/2}$ in the $y$ direction, and each cell of the mesh $I_{i,j}$ is taken as $ [x_{i-1/2},x_{i+1/2}] \times [y_{j-1/2},y_{j+1/2}]$ with its cell center $(x_i,y_j)=(\frac{x_{i-1/2}+x_{i+1/2}}{2}, \frac {y_{j-1/2} +y_{j+1/2}}{2})$. In the next procedures, $x_i+a\Delta x$ is defined as $x_{i+a}$, while $y_j+b\Delta y$  is set  as $y_{j+b}$.

To design a HWENO scheme, we multiply the equation (\ref{EQ2}) by $\frac{1}{\Delta x\Delta y}$,  $\frac {x-x_i} {(\Delta x)^2\Delta y}$ and $\frac {y-y_j} {(\Delta y)^2\Delta x}$ on both sides, respectively, and we integrate them over $I_{i,j}$, then, apply the integration by parts and employ the numerical flux to approximate the values of the flux at the points on the interface of $I_{i,j}$, lastly, we get the semi-discrete finite volume HWENO scheme, and the  explicit formulas are given as follows,
\begin{equation}
\label{ode2}
\left\{
\begin{aligned}
\frac{d \overline u_{i,j}(t)}{dt}&=-\frac{1} {\Delta x \Delta y}  \int_{y_j-1/2}^{y_j+1/2}[\hat f(u(x_{i+1/2},y))-\hat f(u(x_{i-1/2},y))]dy\\
&-\frac{1} {\Delta x \Delta y}  \int_{x_i-1/2}^{x_i+1/2}[\hat g(u(x,y_{j+1/2}))-\hat g(u(x,y_{j-1/2}))]dx,\\
\frac{d \overline v_{i,j}(t)}{dt}&=-\frac{1} {2\Delta x \Delta y}  \int_{y_j-1/2}^{y_j+1/2}[\hat f(u(x_{i-1/2},y))+\hat f(u(x_{i+1/2},y))]dy+\frac 1{{\Delta x}^2 \Delta y}\int_{I_{i,j}}f(u)dxdy\\
&-\frac{1} {\Delta x \Delta y}  \int_{x_i-1/2}^{x_i+1/2}\frac{(x-x_i)}{\Delta x}[\hat g(u(x,y_{j+1/2}))-\hat g(u(x,y_{j-1/2}))]dx,\\
\frac{d \overline w_{i,j}(t)}{dt}&=-\frac{1} {\Delta x \Delta y}  \int_{y_j-1/2}^{y_j+1/2}\frac{(y-y_j)}{\Delta y}[\hat f(u(x_{i+1/2},y))-\hat f(u(x_{i-1/2},y))]dy\\
&-\frac{1} {2\Delta x \Delta y}  \int_{x_i-1/2}^{x_i+1/2}[\hat g(u(x,y_{j-1/2}))+\hat g(u(x,y_{j+1/2}))]dx+\frac 1{{\Delta x \Delta y}^2}\int_{I_{i,j}}g(u)dxdy,
\end{aligned}
\right.
\end{equation}
with initial conditions $ \overline u_{i,j}(0) =$$ \frac 1 {\Delta x \Delta y} \int_{I_{i,j}} u_0(x,y) dxdy$, $ \overline v_{i,j}(0) =$$ \frac 1 {\Delta x \Delta y} \int_{I_{i,j}} u_0(x,y) \frac{x-x_i} { \Delta x} dxdy$ and $ \overline w_{i,j}(0) =$$ \frac 1 {\Delta x \Delta y} \int_{I_{i,j}} u_0(x,y) \frac{y-y_j} { \Delta y} dxdy$. $\overline u_{i,j}(t)$ is the cell average set as $\frac 1 {\Delta x \Delta y} $$\int_{I_{i,j}} u(x,y,t)dxdy$; $\overline v_{i,j}(t)$ is the first order moment in the  $x$ direction defined as $\frac 1 {\Delta x \Delta y} $$\int_{I_{i,j}} u(x,y,t) \frac{x-x_i} { \Delta x} dxdy$ and $\overline w_{i,j}(t)$ is the first  moment in the  $y$ direction taken as $\frac 1 {\Delta x \Delta y} $$\int_{I_{i,j}} u(x,y,t) \frac{y-y_j} { \Delta y} dxdy$. $\hat f(u(x_{i+1/2},y))$ is a numerical flux which is the approximation to the values of  the numerical flux $ f(u) $ at the interface point $(x_{i+1/2},y)$ and $ \hat g(u(x,y_{j+1/2}))$ is a numerical flux to approximate the values of $g(u)$ at  the interface point $(x,y_{j+1/2})$.

Just as in the one dimensional case, we will approximate the integral terms of equations (\ref{ode2}) by employing numerical integration. Since we construct a fourth-order accuracy scheme, 2-point Gaussian will be used in each numerical quadrature, then, these approximated formulas for the integral terms are given as follows,
\begin{equation*}
 \frac 1{\Delta x\Delta y} \int_{I_{i,j}}f(u)dxdy \approx \sum_{k=1}^{2}\sum_{l=1}^{2}\omega_k\omega_l f(u(x_{G_k},y_{G_l})),
\end{equation*}
\begin{equation*}
  \int_{y_j-1/2}^{y_j+1/2}\hat f(u(x_{i+1/2},y))dy\approx \Delta y \sum_{k=1}^{2} \omega_k \hat f(u(x_{i+1/2},y_{G_k})),
\end{equation*}
\begin{equation*}
  \int_{x_i-1/2}^{x_i+1/2}\frac{(x-x_i)}{\Delta x}\hat g(u(x,y_{j+1/2}))dx\approx \Delta x \sum_{k=1}^{2} \omega_k  \frac{(x_{G_k}-x_i)}{\Delta x}\hat g(u(x_{G_k},y_{j+1/2})),
\end{equation*}
and the approximated expressions  for other integral terms are similar. Here, $\omega_1=\frac12$ and $\omega_2=\frac12$ are the quadrature weights, and the coordinates of the Gaussian points over the cell $I_{i,j}$ are
\begin{equation*}
x_{G_1}=x_{i-\sqrt 3 /6},\quad x_{G_2}=x_{i+\sqrt 3 /6},\quad y_{G_1}=y_{j-\sqrt 3 /6},\quad y_{G_2}=y_{j+\sqrt 3 /6}.
\end{equation*}

In the two dimensional case, two fluxes in the $x$ direction and in the $y$ direction are approximated by Lax-Friedrichs numerical flux:
\begin{equation*}\label{2dflx}
  \hat f(u(G_b))=\frac 1 2[f(u^-(G_b))+f(u^+(G_b))]-\frac{\alpha}2(u^+(G_b)-u^-(G_b)),
\end{equation*}
and
\begin{equation*}\label{2dfly}
\hat g(u(G_b))=\frac 1 2[g(u^-(G_b))+g(u^+(G_b))]-\frac{\beta}2(u^+(G_b)-u^-(G_b)),
\end{equation*}
where $\alpha= \max_u|f'(u)|$, $\beta= \max_u|g'(u)|$, and $G_b$ is the Gaussian point on the interface of the cell $I_{i,j}$.

 The general frameworks for the hybrid HWENO scheme are:
 \Red{in Steps 4 and 5, we present the spatial reconstruction for the semi-discrete scheme  (\ref{ode2}). In Step 6, the equations (\ref{ode2}) is discretized by the third order TVD Runge-Kutta method \cite{so1} in time.}

\textbf{Step 4.} Identify the troubled-cell  and  modify the first order moments in the troubled-cell.

We also use the KXRCF troubled-cell indicator by Krivodonova et al. \cite{LJJN} (KXRCF) to identify the discontinuities, which has been introduced in the one dimensional problems, and its explicit expression can be seen in (\ref{indicator}). In particular, the troubled-cell indicator works separately in the $x$ and $y$ directions for two dimensional case, and the norm is still $L_\infty$ norm in the two dimensional case. The cell $I_{i,j}$ is finally identified  as a troubled-cell, if it is identified either in $x$ direction or $y$ direction. In addition, if the cell $I_{i,j}$ is identified as a troubled cell, we  mark the cells $I_{i-1,j}$, $I_{i+1,j}$, $I_{i,j-1}$ and $I_{i,j+1}$  as  troubled cells in practice, as the spatial reconstruction for the neighbor cells also need to use the information of $I_{i,j}$.

If the cell $I_{i,j}$ is identified as a troubled cell, we would modify the first order moments $\overline v_{i,j}$ and $\overline w_{i,j}$ following as next procedures. We modify the first order moments in the troubled cells using dimensional by dimensional manner, and the modification procedures are the same as the one dimensional case. More explicitly, if the cell $I_{i,j}$ is identified as a troubled cell, we use these information $\overline u_{i-1,j}$, $\overline u_{i,j}$, $\overline u_{i+1,j}$, $\overline v_{i-1,j}$, $\overline v_{i+1,j}$ to  reconstruct  the value of  $\overline v_{i,j}$, and the procedures are the same as the expressions introduced in Step 1 for one dimensional case, and the procedures for the  modification of $\overline w_{i,j}$ are similar.

\textbf{Step 5.} Reconstruct the point values of the solutions $u$ at the specific points.

 This step is to reconstruct the point values of $u^+(x_{i-1/2},y_{j\pm\sqrt 3/6})$, $u^-(x_{i+1/2},y_{j\pm\sqrt 3/6})$, $u^+(x_{i\pm\sqrt 3/6},y_{j-1/2})$, $u^-(x_{i\pm\sqrt 3/6},y_{j+1/2})$
and $u(x_{i\pm\sqrt 3/6},y_{j\pm\sqrt 3/6})$ in the cell $I_{i,j}$.
If the cell is identified as a troubled cell, in   Step 5.1, the interface points values of the cell $I_{i,j}$ are reconstructed by HWENO methodology but the internal points values of $I_{i,j}$ are approximated by linear approximation, respectively; otherwise, we directly use the  linear approximation presented in  Step 5.2.

\textbf{Step 5.1.} Reconstruct the point values of the solutions $u$ at the interface points by HWENO approximation  and approximate the internal points  values using linear approximation.

\begin{figure}
  \centering
   \begin{tikzpicture}
     \draw(0,0)rectangle+(1.5,1.5);\draw(1.5,0)rectangle+(1.5,1.5);\draw(3,0)rectangle+(1.5,1.5);
     \draw(0,1.5)rectangle+(1.5,1.5);\draw(1.5,1.5)rectangle+(1.5,1.5);\draw(3,1.5)rectangle+(1.5,1.5);
     \draw(0,3)rectangle+(1.5,1.5);\draw(1.5,3)rectangle+(1.5,1.5);\draw(3,3)rectangle+(1.5,1.5);
     \draw(0.75,0.75)node{1};\draw(2.25,0.75)node{2};\draw(3.75,0.75)node{3};
     \draw(0.75,2.25)node{4};\draw(2.25,2.25)node{5};\draw(3.75,2.25)node{6};
     \draw(0.75,3.75)node{7};\draw(2.25,3.75)node{8};\draw(3.75,3.75)node{9};
     \draw(0.75,-0.25)node{i-1};\draw(2.25,-0.25)node{i};\draw(3.75,-0.25)node{i+1};
     \draw(4.9,0.75)node{j-1};\draw(4.9,2.25)node{j};\draw(4.9,3.75)node{j+1};
     \draw(1.5,2.25-1.732/6.0*1.5)node{*};\draw(1.5,2.25+1.732/6.0*1.5-0.08)node{*};
     \draw(3.0,2.25-1.732/6.0*1.5)node{*};\draw(3.0,2.25+1.732/6.0*1.5-0.08)node{*};
     \draw(2.25-1.732/6.0*1.5,1.45)node{*};\draw(2.25+1.732/6.0*1.5,1.45)node{*};
     \draw(2.25-1.732/6.0*1.5,2.95)node{*};\draw(2.25+1.732/6.0*1.5,2.95)node{*};
     \draw(2.25-1.732/6.0*1.5,2.25-1.732/6.0*1.5)node{*};\draw(2.25+1.732/6.0*1.5,2.25-1.732/6.0*1.5)node{*};
     \draw(2.25-1.732/6.0*1.5,2.25+1.732/6.0*1.5-0.08)node{*};\draw(2.25+1.732/6.0*1.5,2.25+1.732/6.0*1.5-0.08)node{*};
   \end{tikzpicture}
 \caption{The big stencil and its new labels.}
 \label{2dbig}
\end{figure}
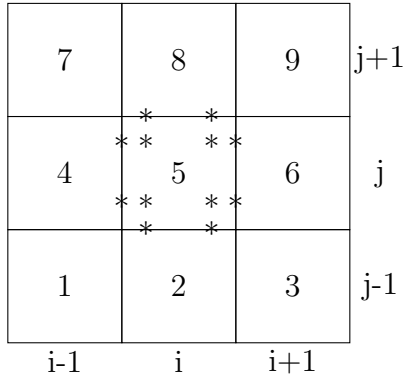

\begin{figure}
  \centering
   \begin{tikzpicture}
   \draw(0+0,0+14.4)       rectangle+(1.8,1.8);\draw(1.8+0,0+14.4)     rectangle+(1.8,1.8);
     \draw(0+0,1.8+14.4)   rectangle+(1.8,1.8);\draw(1.8+0,1.8+14.4)   rectangle+(1.8,1.8);
    \draw(0.9+0,0.9+14.4)node{1};\draw(2.7+0,0.9+14.4)node{2};
    \draw(0.9+0,2.7+14.4)node{4};\draw(2.7+0,2.7+14.4)node{5};
     \draw(0.9+0,-0.25+14.4)node{i-1};\draw(2.7+0,-0.25+14.4)node{i};
    \draw(4.0+0,0.9+14.4)node{j-1};\draw(4.0+0,2.7+14.4)node{j};
   \draw(0+5.0,0+14.4) rectangle+(1.8,1.8);\draw(1.8+5.0,0+14.4)rectangle+(1.8,1.8);
   \draw(0+5.0,1.8+14.4)rectangle+(1.8,1.8);\draw(1.8+5.0,1.8+14.4)rectangle+(1.8,1.8);
   \draw(0.9+5.0,0.9+14.4)node{2};\draw(2.7+5.0,0.9+14.4)node{3};
   \draw(0.9+5.0,2.7+14.4)node{5};\draw(2.7+5.0,2.7+14.4)node{6};
   \draw(0.9+5.0,-0.25+14.4)node{i};\draw(2.7+5.0,-0.25+14.4)node{i+1};
   \draw(4.0+5.0,0.9+14.4)node{j-1};\draw(4.0+5.0,2.7+14.4)node{j};
     \draw(0+0,0+9.6)       rectangle+(1.8,1.8);
     \draw(1.8+0,0+9.6)     rectangle+(1.8,1.8);
     \draw(0+0,1.8+9.6)   rectangle+(1.8,1.8);
     \draw(1.8+0,1.8+9.6)   rectangle+(1.8,1.8);
    \draw(0.9+0,0.9+9.6)node{4};\draw(2.7+0,0.9+9.6)node{5};
  \draw(0.9+0,2.7+9.6)node{7}; \draw(2.7+0,2.7+9.6)node{8};
     \draw(0.9+0,-0.25+9.6)node{i-1};\draw(2.7+0,-0.25+9.6)node{i};
    \draw(4+0,0.9+9.6)node{j};\draw(4+0,2.7+9.6)node{j+1};
   \draw(0+5.0,0+9.6)       rectangle+(1.8,1.8);\draw(1.8+5.0,0+9.6)     rectangle+(1.8,1.8);
     \draw(0+5.0,1.8+9.6)   rectangle+(1.8,1.8);\draw(1.8+5.0,1.8+9.6)   rectangle+(1.8,1.8);
    \draw(0.9+5.0,0.9+9.6)node{5}; \draw(2.7+5.0,0.9+9.6)node{6};
    \draw(0.9+5.0,2.7+9.6)node{8};\draw(2.7+5.0,2.7+9.6)node{9};
     \draw(0.9+5.0,-0.25+9.6)node{i};\draw(2.7+5.0,-0.25+9.6)node{i+1};
    \draw(4+5.0,0.9+9.6)node{j};\draw(4+5.0,2.7+9.6)node{j+1};
     \draw(0.0+0,0+4.8)rectangle+(1.2,1.2);\draw(1.2+0,0+4.8)rectangle+(1.2,1.2);\draw(2.4+0,0+4.8)rectangle+(1.2,1.2);
     \draw(0+0,1.2+4.8)rectangle+(1.2,1.2);\draw(1.2+0,1.2+4.8)rectangle+(1.2,1.2);\draw(0+0,2.4+4.8)rectangle+(1.2,1.2);
     \draw(0.6+0,0.6+4.8)node{1};\draw(1.8+0,0.6+4.8)node{2};\draw(3+0,0.6+4.8)node{3};
     \draw(0.6+0,1.8+4.8)node{4};\draw(1.8+0,1.8+4.8)node{5};\draw(0.6+0,3+4.8)node{7};
     \draw(0.6+0,-0.25+4.8)node{i-1};\draw(1.8+0,-0.25+4.8)node{i};\draw(3.0+0,-0.25+4.8)node{i+1};
     \draw(4.0+0,0.6+4.8)node{j-1};\draw(4.0+0,1.8+4.8)node{j};\draw(4.0+0,3.0+4.8)node{j+1};
     \draw(0.0+5,0+4.8)rectangle+(1.2,1.2);\draw(1.2+5,0+4.8)rectangle+(1.2,1.2);\draw(2.4+5,0+4.8)rectangle+(1.2,1.2);
     \draw(1.2+5,1.2+4.8)rectangle+(1.2,1.2);\draw(2.4+5,1.2+4.8)rectangle+(1.2,1.2);\draw(2.4+5,2.4+4.8)rectangle+(1.2,1.2);
     \draw(0.6+5,0.6+4.8)node{1};\draw(1.8+5,0.6+4.8)node{2};\draw(3+5,0.6+4.8)node{3};
     \draw(1.8+5,1.8+4.8)node{5};\draw(3+5,1.8+4.8)node{6};\draw(3+5,3+4.8)node{9};
     \draw(0.6+5,-0.25+4.8)node{i-1};\draw(1.8+5,-0.25+4.8)node{i};\draw(3.0+5,-0.25+4.8)node{i+1};
     \draw(4.0+5,0.6+4.8)node{j-1};\draw(4.0+5,1.8+4.8)node{j};\draw(4.0+5,3.0+4.8)node{j+1};
     \draw(0.0+0,0+0)rectangle+(1.2,1.2);\draw(0+0,1.2+0)rectangle+(1.2,1.2);\draw(1.2+0,1.2+0)rectangle+(1.2,1.2);
     \draw(0+0,2.4+0)rectangle+(1.2,1.2);\draw(1.2+0,2.4+0)rectangle+(1.2,1.2);\draw(2.4+0,2.4+0)rectangle+(1.2,1.2);
     \draw(0.6+0,0.6+0)node{1};\draw(0.6+0,1.8+0)node{4};\draw(1.8+0,1.8+0)node{5};
     \draw(0.6+0,3+0)node{7};\draw(1.8+0,3+0)node{8};\draw(3+0,3+0)node{9};
     \draw(0.6+0,-0.25+0)node{i-1};\draw(1.8+0,-0.25+0)node{i};\draw(3.0+0,-0.25+0)node{i+1};
     \draw(4+0,0.6+0)node{j-1};\draw(4+0,1.8+0)node{j};\draw(4+0,3.0+0)node{j+1};
     \draw(2.4+5,0+0)rectangle+(1.2,1.2);\draw(1.2+5,1.2+0)rectangle+(1.2,1.2);
     \draw(2.4+5,1.2+0)rectangle+(1.2,1.2);\draw(0+5,2.4+0)     rectangle+(1.2,1.2);
     \draw(1.2+5,2.4+0)rectangle+(1.2,1.2);\draw(2.4+5,2.4+0)rectangle+(1.2,1.2);
     \draw(3+5,0.6+0)node{3};\draw(1.8+5,1.8+0)node{5};\draw(3+5,1.8+0)node{6};
     \draw(0.6+5,3+0)node{7};\draw(1.8+5,3+0)node{8};\draw(3+5,3+0)node{9};
     \draw(0.6+5,-0.25+0)node{i-1};\draw(1.8+5,-0.25+0)node{i};\draw(3.0+5,-0.25+0)node{i+1};
     \draw(4+5,0.6+0)node{j-1};\draw(4+5,1.8+0)node{j};\draw(4+5,3.0+0)node{j+1};
   \end{tikzpicture}
 \caption{The eight small stencils and these respective labels. From left to right and top to bottom are the stencils: $S_1,...,S_8$.}
 \label{2dsmall}
\end{figure}
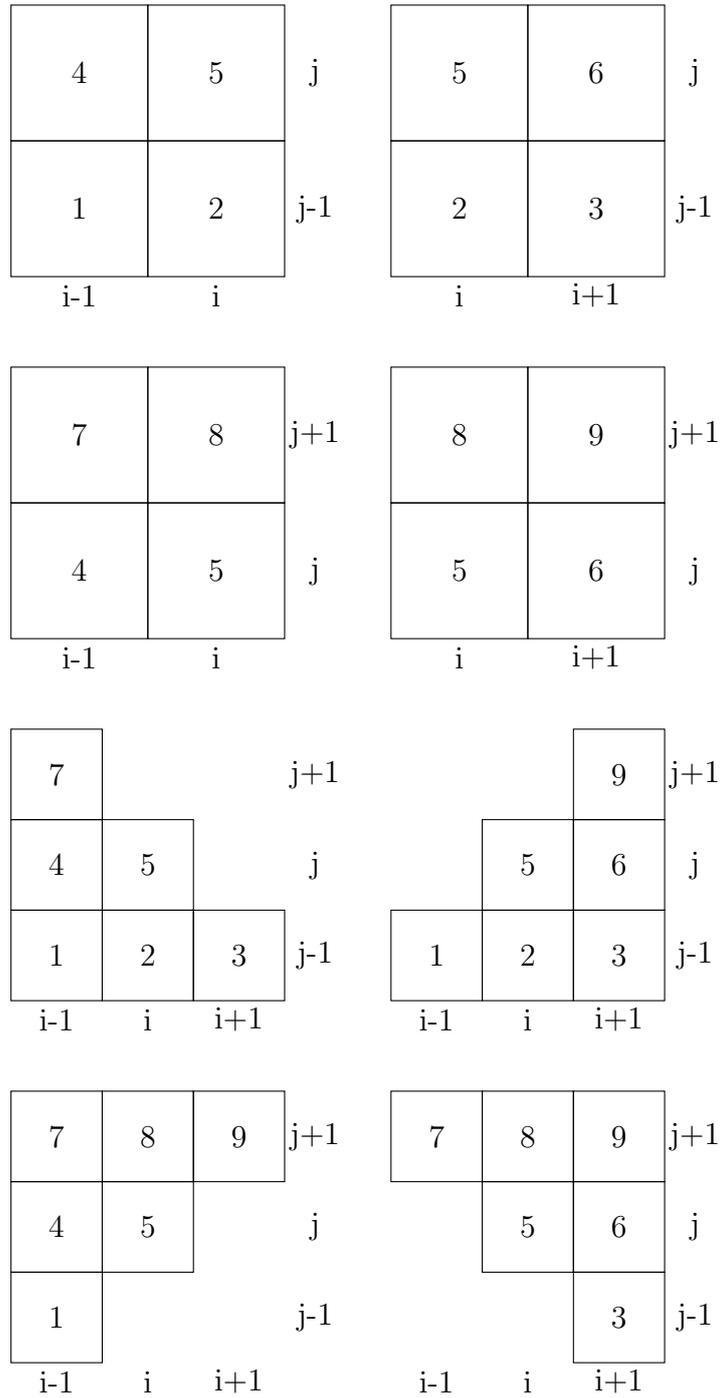

We first give the big stencil $S_0$ in Figure \ref{2dbig}. For simplicity, we rebel  the cell $I_{i,j}$ and its neighboring cells  as $I_1,...,I_9$. Particularly, the new label of the cell $I_{i,j}$ is $I_5$ and the symbols "*" in $I_5$ represent the locations of the solutions $u$ where we need to  reconstruct. We also give eight small stencils: $S_1,S_2,...,S_8$ shown in Figure \ref{2dsmall}, which were first introduced by Qiu and Shu in \cite{QSHw}. Noticed that we reconstruct the point values of solutions $u$ in the cell $I_{i,j}$, then, we would like to use more information in the cell $I_{i,j}$, such as the first order moments $ \overline v_{i,j}$ and $ \overline w_{i,j}$, so we construct eight incomplete cubic reconstruction polynomials, and these polynomials have the expressions as follows
\begin{equation}\label{2drepoly}
  p_n(x,y)= a_0+a_1x+a_2y+a_3x^2+a_4xy+a_5y^2+a_6x^3+a_7y^3,\quad n=1,...,8,
\end{equation}
satisfying as
  \begin{equation*}
\begin{array}{ll}
\frac{1}{\Delta x \Delta y}\int_{I_k}p_n(x,y)dxdy=\overline u_k,            \\
 \frac{1}{\Delta x \Delta y}\int_{I_{k_x}}p_n(x,y)\frac{(x-x_{k_x})}{\Delta x}dxdy=\overline v_{k_x}, \quad \frac{1}{\Delta x \Delta y}\int_{I_{k_y}}p_n(x,y)\frac{(y-y_{k_y})}{\Delta y}dxdy=\overline w_{k_y}, \\
\end{array}
\end{equation*}
for
  \begin{equation*}
  \begin{array}{ll}
n=1\ k=1,2,4,5,\ k_x=4,5,\ k_y=2,5; \quad n=2,\ k=2,3,5,6,\ k_x=5,6,\ k_y=2,5;   \\
n=3\ k=4,5,7,8,\ k_x=4,5,\ k_y=5,8; \quad n=4,\ k=5,6,8,9,\ k_x=5,6,\ k_y=5,8;   \\
n=5\ k=1,2,3,4,5,7,\ k_x=5,\ k_y=5; \quad n=6,\ k=1,2,3,5,6,9,\ k_x=5,\ k_y=5;   \\
n=7\ k=1,4,5,7,8,9,\ k_x=5,\ k_y=5; \quad n=8,\ k=3,5,6,7,8,9,\ k_x=5,\ k_y=5.
\end{array}
\end{equation*}
Then, we combine these eight incomplete cubic polynomials to obtain a fourth-order approximation for the reconstruction  of these points values of the solutions $u$. For simplicity, we use $G_k$ to represent the specific points, where we want to reconstruct. At first, we use the linear weights $\gamma_1^{k},...,\gamma_8^{k}$  to joint these eight small polynomials, satisfying as
 \begin{equation}\label{2dpointr}
u(G_k)=\sum_{n=1}^{8} \gamma_n^k p_n(G_k).
\end{equation}
If $\sum_{n=1}^{8}\gamma_n^{k}=1$, the requirement (\ref{2dpointr}) are always satisfied for any incomplete cubic polynomial $u$, and the form of the polynomial is presented in (\ref{2drepoly}), but we still have two other constraints on the linear weights to hold the requirement (\ref{2dpointr}) for $u=x^2y$ and $xy^2$. Subject to these three constraints listed above, it leaves five free parameters to calculate the linear weights, and we can obtain $\gamma_1^{k},...,\gamma_8^{k}$ easily and  uniquely  by  minimizing $\sum_{n=1}^{8} (\gamma_n^{k})^2$. In fact, the linear weights $\gamma_1^{k},...,\gamma_8^{k}$ determined  by this least square methodology  are all positive in the implementation. For simplicity, we only present the eight linear weights to reconstruct $u^-(x_{i+1/2},y_{j+\sqrt 3/6})$ at the interface point, and the values are $\frac{3533+351\sqrt 3}{37040}$, $\frac{5727+351\sqrt 3}{37040}$, $\frac{3533-351\sqrt 3}{37040}$, $\frac{5727-351\sqrt 3}{37040}$,
$\frac{10599-1867\sqrt 3}{111120}$, $\frac{17181-415\sqrt 3}{111120}$, $\frac{10599+1867\sqrt3}{111120}$, $\frac{17181+415\sqrt 3}{111120}$, respectively, and the linear weights for other points on the interface can be obtained by symmetry. In addition, it is interesting that the linear weights to reconstruct $u(x_{i\pm\sqrt 3/6},y_{j\pm\sqrt 3/6})$ are all $\frac 18$.

Similarly as in the one dimensional problems, if $G_k$ is inside of $I_{i,j}$, we
directly use linear approximation to reconstruct $u(G_k)$ as $\sum_{n=1}^{8} \gamma_n^{k}p_n(G_k)$, and we'd better  simplify the formals  in advance, instead of calculating in the codes; while $G_k$ is located on the interface of the cell $I_{i,j}$, we need to employ the next HWENO reconstruction procedures, then, we first compute the smoothness indicators $\beta_n$, which measure how smooth the function $p_n(x,y)$ in the cell $I_{i,j}$. The formula was listed by \cite{hs}, given as follows,
\begin{equation}\label{2dGHYZ}
  \beta_n= \sum_{|l|=1}^3|I_{i,j}|^{|l|-1} \int_{I_{i,j}}\left( \frac {\partial^{|l|}}{\partial x^{l_1}\partial y^{l_2}}p_n(x,y)\right)^2 dxdy, \quad n=1,...8,
\end{equation}
where $l=(l_1,l_2)$, $|l|=l_1+l_2$, then, we can get the non-linear weights using the linear weights and  the smoothness indicators, having
\begin{equation*}
\label{2d91}
\omega_n^{k}=\frac{\bar\omega_n^{k}}{\sum_{l=1}^{8}\bar\omega_{l}^{k}},
\ \mbox{with} \ \bar\omega_{n}^{k}=\frac {\gamma_{n}^{k}}{(\beta_{n}+\varepsilon)^2}, \ n=1,...,8,
\end{equation*}
where $\varepsilon$ is set as $10^{-6}$ just as in one dimensional case. Finally, the approximation for the point values of the solutions $u$ at the interface  point $G_k$ is given by
\begin{equation*}
  u^*(G_k) = \sum_{n=1}^{8} \omega_n^{k}p_n(G_k),
\end{equation*}
where "*" is "+" when $G_k$ is located on the left or bottom interface of the cell $I_{i,j}$, while "*" is "-" on the right or top interface of $ I_{i,j} $.

\textbf{Step 5.2.} Reconstruct the point values of the solutions $u$ at the specific points by linear approximation straightforwardly.

In this step, we'll use the same polynomials and linear weights introduced in  Step 5.1, then, the linear approximation of the solutions $u$ at reconstructed point $G_k$ can be directly taken  as
\begin{equation*}
  u^*(G_k) = \sum_{n=1}^{8} \gamma_n^{(k)}p_n(G_k).
\end{equation*}
If $G_k$ is located on the interface of the cell $I_{i,j}$, "*" has the same meaning just as in Step 5.1; otherwise, "*" will be blank. Similarly,  we also can obtain the simplified formulas  easily  for the linear approximation of $u^*(G_k)$ in advance, instead of calculating over and over again in the codes.

 \Red{Noticed that the reconstruction for the points values of the solutions $u$ only has the fourth order accuracy for the information we used here is not enough to reconstruct interpolation polynomial with degree 4, therefore, the scheme only is the fourth order in the two dimensional case.}

\textbf{Step 6.} Discretize the  semi-discrete scheme (\ref{ode2}) in time by the third order TVD Runge-Kutta method \cite{so1}.

When we have finished Steps 4 and 5, the  semi-discrete scheme (\ref{ode2}) is discretized in time by the third order TVD Runge-Kutta method, and the explicit expression has been presented in (\ref{RK}) for the one dimensional case.

{\bf \em Remark 3:} The KXRCF indicator is satisfying for two dimensional hyperbolic conservation laws. For two dimensional scalar equation, the solution $u$ is set as the indicator variable. $\overrightarrow{v}$ is taken as $f'(u)$ in the $x$ direction, while it is defined as $g'(u)$ in the $y$ direction; for two dimensional Euler equations, the density $\rho$ and the energy $E$ are taken as the indicator variables, respectively. $\overrightarrow{v}$ is defined as the velocity $\mu$ in the $x$ direction of the fluid, while it is set as the  velocity $\nu$  in the $y$ direction of the fluid.

{\bf \em Remark 4:} For the systems, such as the two dimensional compressible Euler equations, similarly as in the one dimension case, we first use the KXRCF indicator to identify the troubled cell in  Step 4. If the cell $I_{i,j}$ is identified as a troubled cell, we modify the first order moments for each component. For the modification, it is different from the one dimensional case for it has two first order moments and two flux functions, so we modify the moments $\overline v_{i,j}$ in the $x$ direction in terms of  the  local characteristic direction provided by $f(u)$, while reconstruct $\overline w_{i,j}$ in the $y$ direction based on  the  local characteristic direction of $g(u)$. For each local characteristic direction, we follow the procedures of Step 4 to reconstruct the first order moments in the troubled cells. For the reconstruction for the point values of the solutions $u$,  all HWENO procedures are performed on  the  local characteristic decompositions, and the linear approximations are based on component by component.

\section{Numerical tests}
\label{sec3}
\setcounter{equation}{0}
\setcounter{figure}{0}
\setcounter{table}{0}

In this section, we perform the numerical results of the   hybrid  HWENO scheme in the one and two dimensional cases, which is outlined in  the previous section. If no otherwise specified, the CFL number is set as 0.6 for one dimensional tests and 0.45 for two dimensional examples.

\subsection{Accuracy tests}

For simplicity, Hybrid HWENO scheme is denoted as the hybrid HWENO scheme introduced in the previous section, while HWENO scheme is represented as that we modify the first order moments for every cell and employ HWENO reconstruction at the interface points for the spatial discretization on the basis of the hybrid HWENO scheme. WENO scheme was listed by Jiang and Shu \cite{js}, while Hybrid WENO scheme using KXRCF troubled-cell indicator was introduced by Li and Qiu \cite{Glj}. Since the four schemes all have fifth order accuracy in one dimensional problems, we will compare their performance in the one dimensional  accuracy tests. For two dimensional smooth tests, as the Hybrid HWENO scheme is based on the finite volume framework, we also make the comparisons with the classical fifth order finite volume WENO scheme narrated in \cite{s2}.

\begin{figure}
 \centerline{
\psfig{file=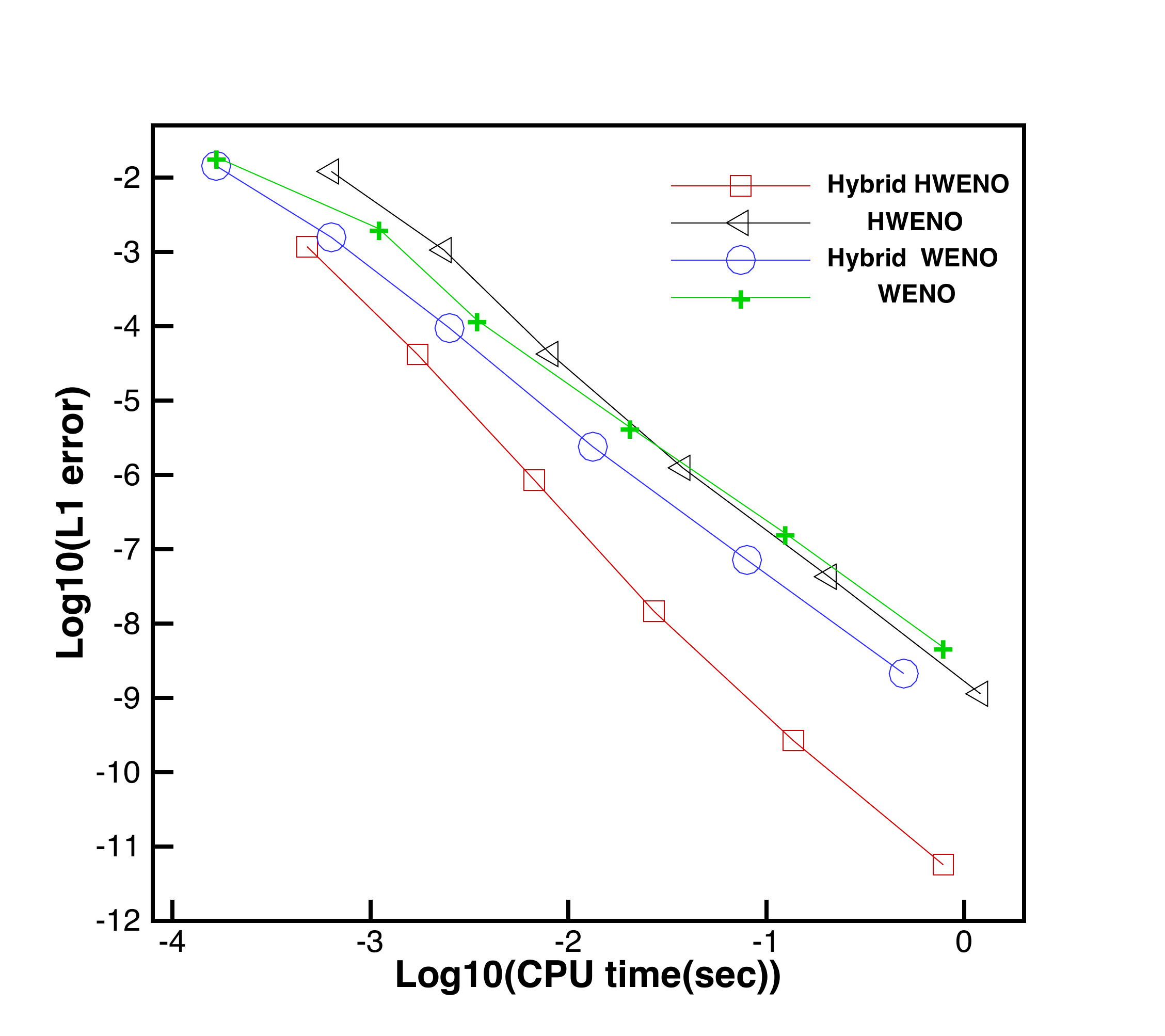,width=2.5 in}
\psfig{file=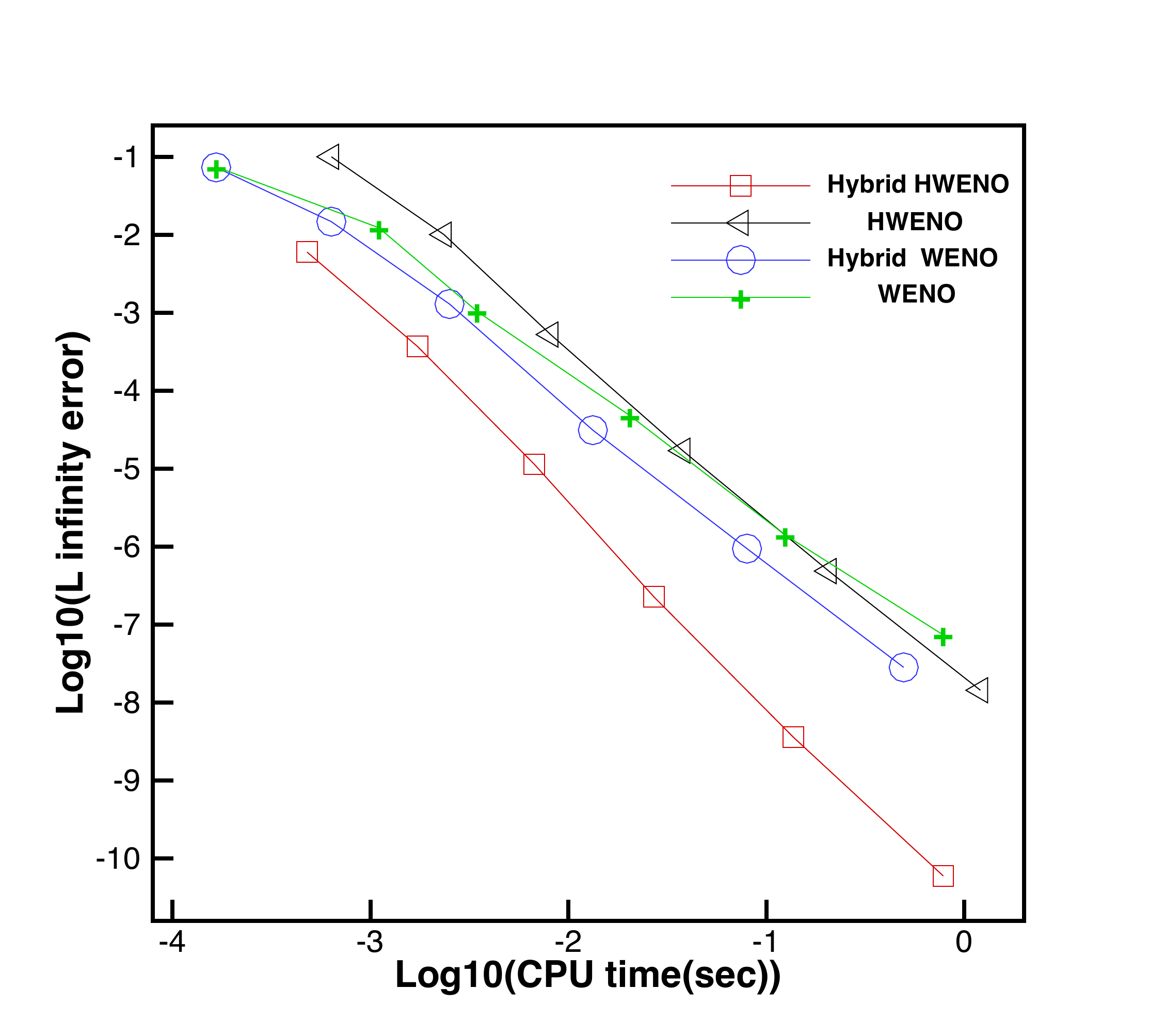,width=2.5 in}}
 \caption{1D-Burgers' equation: initial data
$u(x,0)=0.5+sin(\pi x)$. $T=0.5/\pi$. Computing times and errors. Rectangle signs and a red solid  line denote the results of Hybrid HWENO scheme; triangle signs and a black solid  line denote the results of HWENO scheme; circle signs and a blue solid  line denote the results of  Hybrid WENO scheme; plus signs and a green solid  line denote the results of WENO scheme.}
\label{Fburges1d_smooth}
\end{figure}

\noindent{\bf Example 3.1.} We solve the one dimensional Burgers' equation:
\begin{equation}\label{1dbugers}
  u_t+(\frac {u^2} 2)_x=0, \quad 0<x<2,
\end{equation}
with the initial  condition $u(x,0)=0.5+sin(\pi x)$ and periodic boundary condition. We present the numerical results at $t=0.5/\pi$  when the solution is still smooth, then, the numerical errors and orders are shown in Table \ref{tburgers1d} with $N$ uniform meshes for HWENO and WENO schemes. From the table, we can see that all four schemes have fifth order accuracy. Firstly, we know the hybrid schemes have less errors than the original schemes, meanwhile, we also can find that two HWENO schemes have less errors than corresponding  WENO schemes with  the same number of cells.  In Figure \ref{Fburges1d_smooth}, we show numerical errors against CPU times by using four different schemes, which illustrates Hybrid HWENO scheme has much higher efficiency than other three schemes, meanwhile, the two HWENO schemes only need three cells while the two WENO schemes need five cells for the spatial reconstruction.
\begin{table}
\begin{center}
\caption{1D-Burgers' equation: initial data
$u(x,0)=0.5+sin(\pi x)$. HWENO and WENO schemes. $T=0.5/\pi$. $L^1$ and $L^\infty$ errors and orders. Uniform meshes with $N$ cells. }
\medskip
\begin{tabular} {|c|c|c|c|c|c|c|c|c|} \hline
& \multicolumn{4}{c|}{HWENO  scheme} & \multicolumn{4}{c|}{WENO scheme}\\
\hline\hline
  $N$ cells &$ L^1$ error &  order & $L^\infty$error &  order &$ L^1$ error &  order & $L^\infty$ error &order\\   \hline
   10 &     1.21E-02 &         &     1.00E-01 &         							 &     1.90E-02 &                 &     7.46E-02 &       \\ \hline
     20 &     1.06E-03 &       3.52 &     1.00E-02 &       3.32 							    &     2.06E-03 &       3.20 &     1.23E-02 &       2.60 \\ \hline
     40 &     4.23E-05 &       4.65 &     5.25E-04 &       4.26 							    &     1.22E-04 &       4.08 &     1.05E-03 &       3.55 \\ \hline
     80 &     1.24E-06 &       5.09 &     1.70E-05 &       4.95 							    &     4.36E-06 &       4.80 &     4.78E-05 &       4.46 \\ \hline
    160 &     4.26E-08 &       4.87 &     4.84E-07 &       5.13 							    &     1.64E-07 &       4.74 &     1.41E-06 &       5.09 \\ \hline
    320 &     1.13E-09 &       5.24 &     1.43E-08 &       5.08 							    &     4.78E-09 &       5.10 &     7.35E-08 &       4.26 \\ \hline
\hline
& \multicolumn{4}{c|}{Hybrid HWENO  scheme} & \multicolumn{4}{c|}{Hybrid WENO scheme}\\
\hline\hline
  $N$ cells &$ L^1$ error &  order & $L^\infty$error &  order &$ L^1$ error &  order & $L^\infty$ error &order\\   \hline
  10 &     1.18E-03 &        &     6.00E-03 &        							&     1.44E-02 &        &     7.32E-02 &              \\ \hline
     20 &     4.18E-05 &       4.82 &     3.69E-04 &       4.02 							     &     1.58E-03 &       3.19 &     1.47E-02 &       2.31 \\ \hline
     40 &     8.51E-07 &       5.62 &     1.14E-05 &       5.02 							     &     9.45E-05 &       4.06 &     1.29E-03 &       3.51 \\ \hline
     80 &     1.46E-08 &       5.87 &     2.26E-07 &       5.65 							     &     2.39E-06 &       5.30 &     3.11E-05 &       5.38 \\ \hline
    160 &     2.66E-10 &       5.78 &     3.59E-09 &       5.98 							     &     7.15E-08 &       5.06 &     9.40E-07 &       5.05 \\ \hline
    320 &     5.65E-12 &       5.56 &     5.93E-11 &       5.92 							     &     2.12E-09 &       5.08 &     2.82E-08 &       5.06 \\ \hline
\end{tabular}
\label{tburgers1d}
\end{center}
\end{table}
\smallskip

\begin{figure}
 \centerline{
\psfig{file=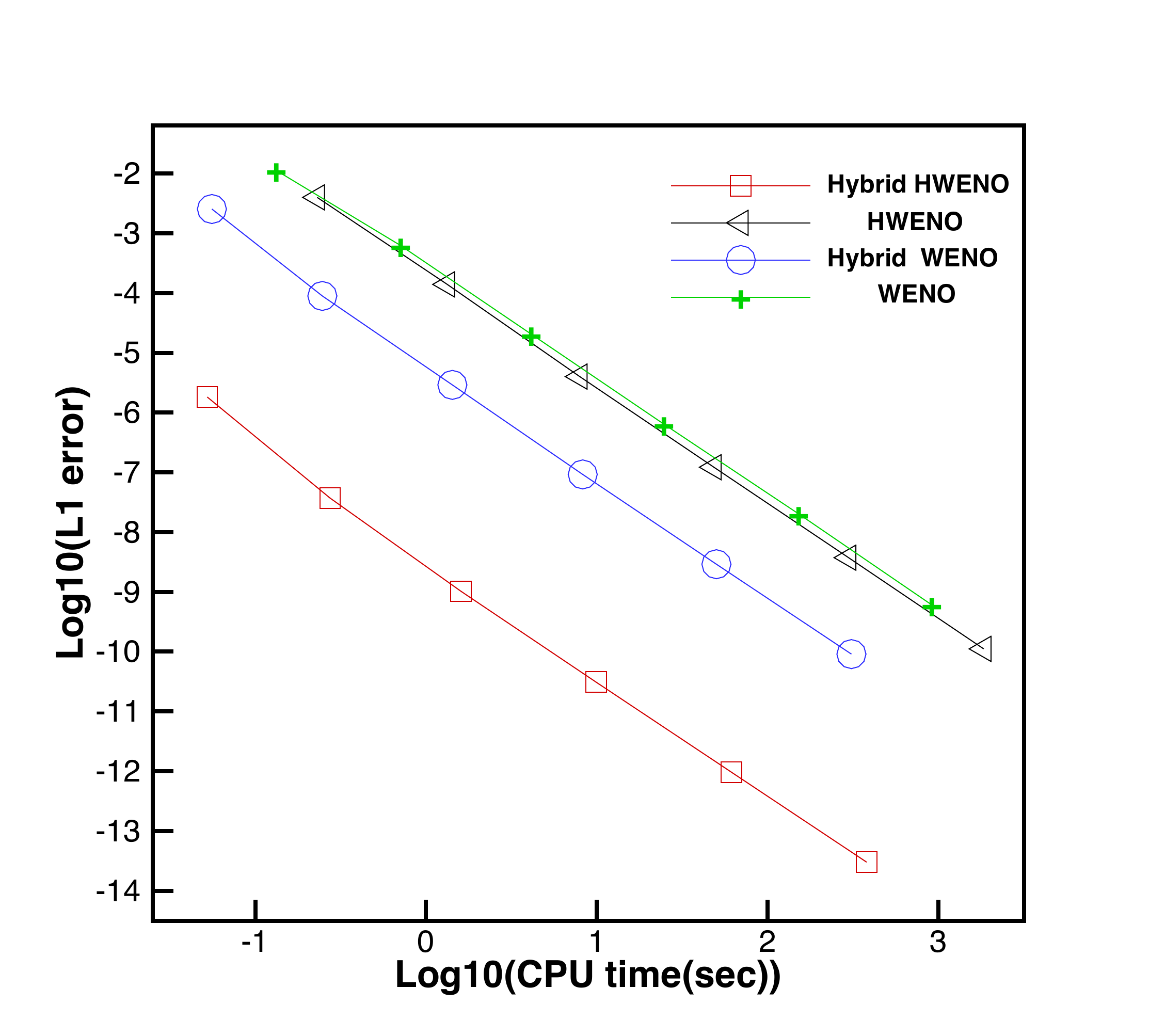,width=2.5 in}
\psfig{file=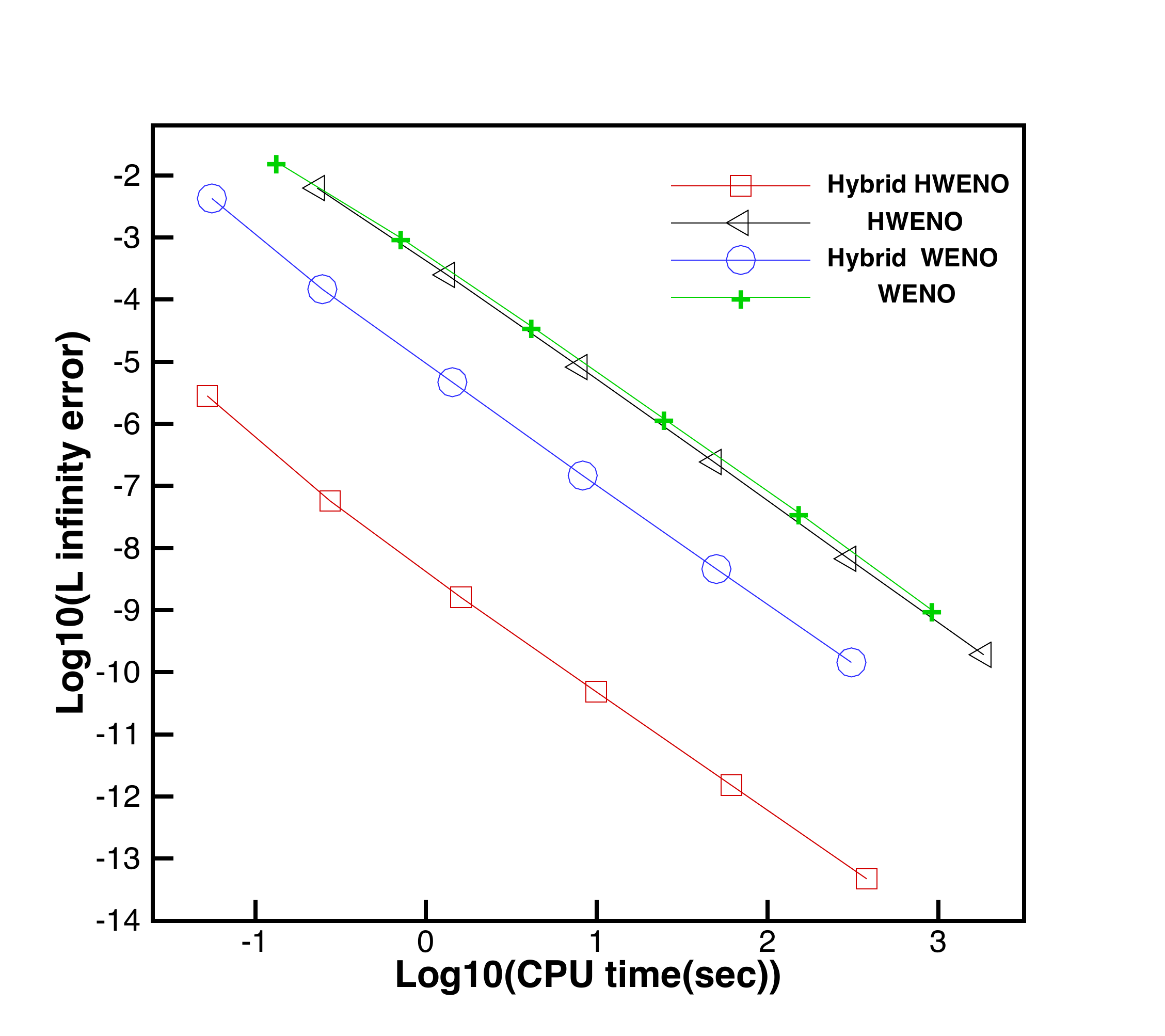,width=2.5 in}}
 \caption{1D-Euler equations: initial data
$\rho(x,0)=1+0.2\sin(\pi x)$, $\mu(x,0)=1$ and $p(x,0)=1$. $T=2$. Computing times and errors. Rectangle signs and a red solid  line denote the results of Hybrid HWENO scheme; triangle signs and a black solid  line denote the results of HWENO scheme; circle signs and a blue solid  line denote the results of  Hybrid WENO scheme; plus signs and a green solid  line denote the results of WENO scheme.}
\label{FEuler1d_smooth}
\end{figure}

\noindent{\bf Example 3.2.} We consider  one dimensional Euler equations:
\begin{equation}
\label{euler1}
 \frac{\partial}{\partial t}
 \left(
 \begin{array}{c}
 \rho \\
 \rho \mu \\
 E
 \end{array} \right )
 +
 \frac{\partial}{\partial x}
 \left (
\begin{array}{c}
\rho \mu \\
\rho \mu^{2}+p \\
\mu(E+p)
\end{array}
 \right )
=0,
\end{equation}
 in which $\rho$ is the density, $\mu$ is the velocity, $E$ is the total energy and $p$ is the pressure. The initial conditions are set to  $\rho(x,0)=1+0.2\sin(\pi x)$, $\mu(x,0)=1$, $p(x,0)=1$ and  $\gamma=1.4$. The computing domain is $ x \in [0, 2\pi]$. Periodic boundary condition is applied here. The exact solution is $\rho(x,t)=1+0.2\sin(\pi(x-t))$, $\mu(x,0)=1$, $p(x,0)=1$ and the final computing time is $T=2$. The numerical errors and  orders of the density for the HWENO and WENO schemes are given in Table \ref{tEluer1d}, which shows four schemes achieve the designed fifth order accuracy. Similarly, the hybrid schemes have less errors than the original schemes and the two HWENO schemes have less errors than corresponding  WENO schemes with  the same number of cells. In addition, Figure \ref{FEuler1d_smooth} represents numerical errors against CPU times using four different schemes, which shows  Hybrid HWENO scheme has  higher efficiency than other three schemes, and the HWENO schemes are more compact than the WENO schemes.
\begin{table}
\begin{center}
\caption{1D-Euler equations: initial data
$\rho(x,0)=1+0.2\sin(\pi x)$, $\mu(x,0)=1$ and $p(x,0)=1$. HWENO and WENO  schemes. $T=2$. $L^1$ and $L^\infty$ errors and orders. Uniform meshes with $N$ cells.}
\medskip
\begin{tabular} {|c|c|c|c|c|c|c|c|c|} \hline
& \multicolumn{4}{c|}{HWENO  scheme} & \multicolumn{4}{c|}{WENO scheme}\\
\hline\hline
  $N$ cells &$ L^1$ error &  order & $L^\infty$error &  order &$ L^1$ error &  order & $L^\infty$ error &order\\   \hline
  10 &     3.98E-03 &        &     6.25E-03  &
        & 1.13E-02 &        &  1.66E-02 &        \\ \hline
     20 &     1.39E-04 &       4.84 &     2.50E-04 &       4.64 								     &     6.26E-04 &       4.17 &     9.94E-04 &       4.06 \\ \hline
     40 &     4.00E-06 &       5.12 &     8.18E-06 &       4.93 							     &     2.04E-05 &       4.94 &     3.72E-05 &       4.74 \\ \hline
     80 &     1.22E-07 &       5.04 &     2.43E-07 &       5.08 								     &     6.45E-07 &       4.98 &     1.21E-06 &       4.94 \\ \hline
    160 &     3.73E-09 &       5.03 &     6.71E-09 &       5.18 								     &     2.01E-08 &       5.01 &     3.67E-08 &       5.05 \\ \hline
    320 &     1.11E-10 &       5.07 &     1.91E-10 &       5.13
    &     6.09E-10 &       5.04 &     1.01E-09 &       5.19 \\ \hline
\hline
& \multicolumn{4}{c|}{Hybrid HWENO  scheme} & \multicolumn{4}{c|}{Hybrid WENO scheme}\\
\hline\hline
  $N$ cells &$ L^1$ error &  order & $L^\infty$error &  order &$ L^1$ error &  order & $L^\infty$ error &order\\   \hline
     10 &     1.82E-06 &            &     2.82E-06 &            							  &     2.55E-03 &            &     4.25E-03 &            \\ \hline
     20 &     3.71E-08 &       5.62 &     5.73E-08 &       5.62 							     &     8.94E-05 &       4.83 &     1.47E-04 &       4.86 \\ \hline
     40 &     1.02E-09 &       5.18 &     1.60E-09 &       5.16 							     &     2.91E-06 &       4.94 &     4.67E-06 &       4.97 \\ \hline
     80 &     3.10E-11 &       5.05 &     4.86E-11 &       5.04 							     &     9.22E-08 &       4.98 &     1.47E-07 &       4.99 \\ \hline
    160 &     9.61E-13 &       5.01 &     1.51E-12 &       5.01 							     &     2.90E-09 &       4.99 &     4.59E-09 &       5.00 \\ \hline
    320 &     3.00E-14 &       5.00 &     4.71E-14 &       5.00 							     &     9.10E-11 &       5.00 &     1.43E-10 &       5.00 \\ \hline
\end{tabular}
\label{tEluer1d}
\end{center}
\end{table}
\smallskip

\begin{figure}
 \centerline{
\psfig{file=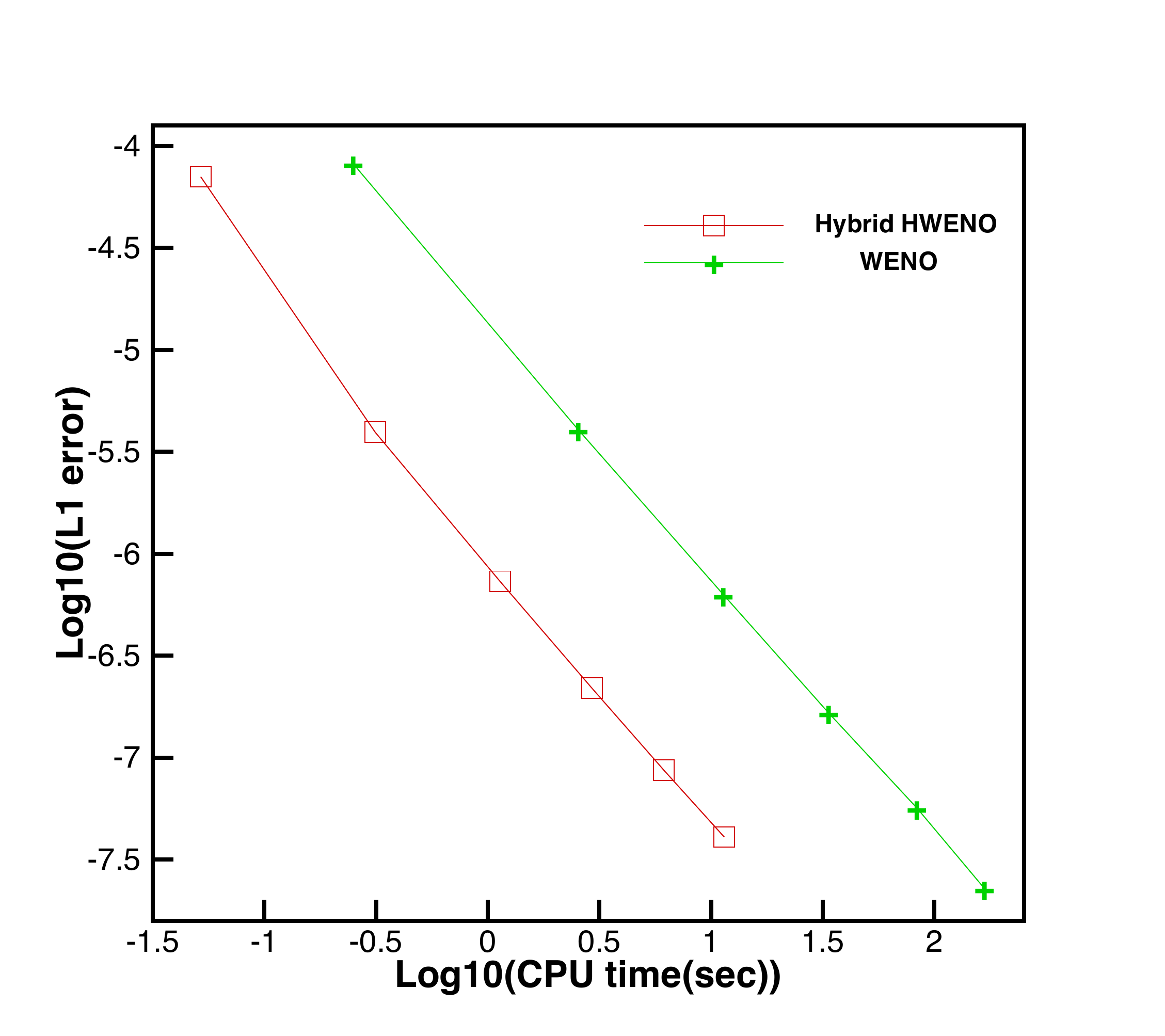,width=2.5 in}
\psfig{file=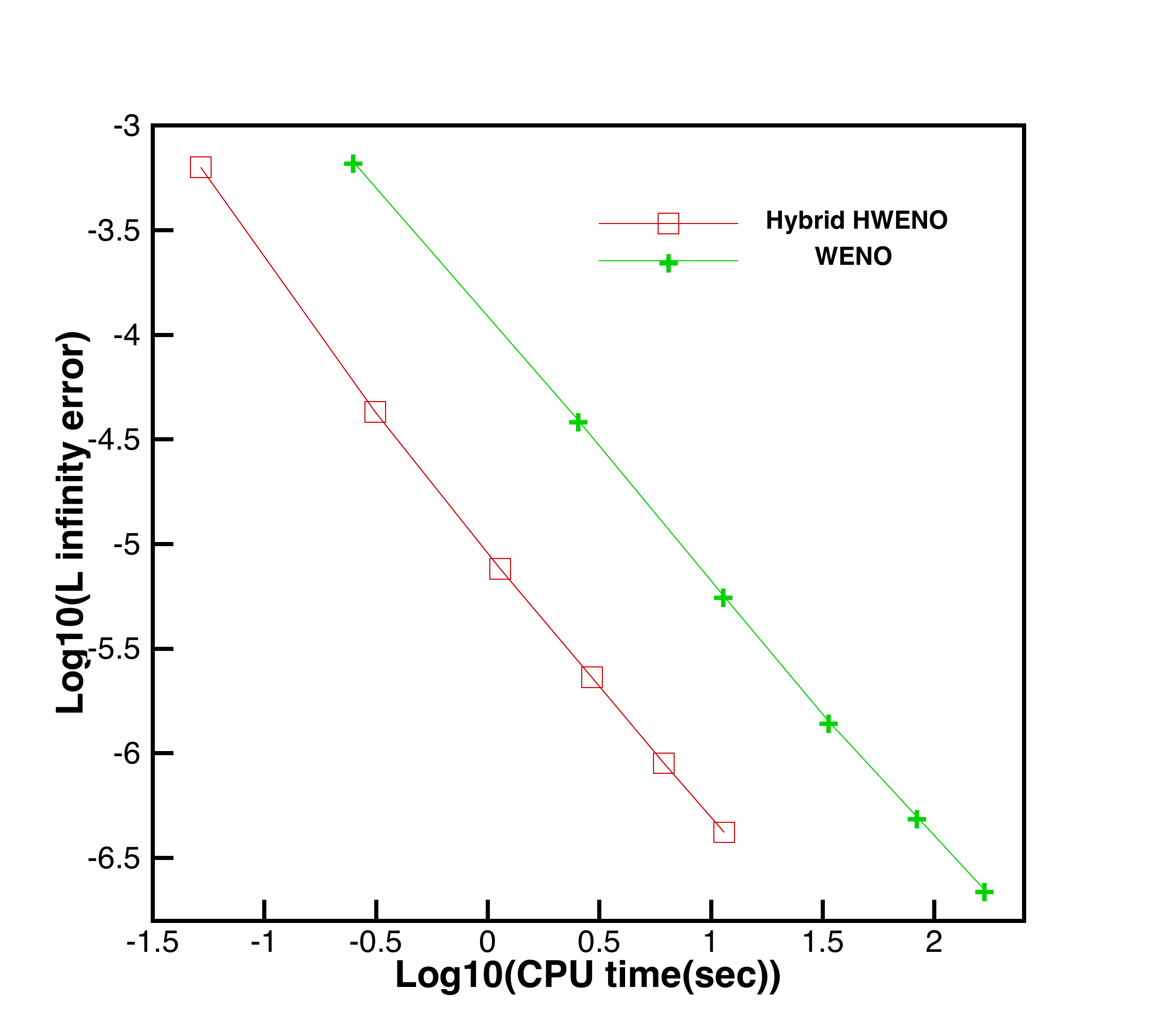,width=2.5 in}}
 \caption{2D-Burgers' equation: initial data
$u(x,y,0)=0.5+sin(\pi (x+y)/2)$. $T=0.5/\pi$. Computing times and errors. Rectangle signs and a red solid  line denote the results of Hybrid HWENO scheme; plus signs and a green solid  line denote the results of WENO scheme.}
\label{Fburges2d_smooth}
\end{figure}

\noindent{\bf Example 3.3.} Consider the two dimensional Burgers' equation:
\begin{equation}\label{2dbugers}
  u_t+(\frac {u^2} 2)_x+(\frac {u^2} 2)_y=0, \quad 0<x<4, \ 0<y<4
\end{equation}
with the initial  condition $u(x,y,0)=0.5+sin(\pi (x+y)/2)$ and periodic boundary condition in each direction. We perform the numerical results at $t=0.5/\pi$. In this time, the solution is still smooth, then, \Red{the numerical errors and order are shown in Table \ref{tburgers2d} for Hybrid HWENO and WENO schemes. We can see that both  schemes achieve the designed order accuracy. In Figure \ref{Fburges2d_smooth}, we present their numerical errors against CPU times, which illustrates Hybrid HWENO scheme has higher efficiency than WENO scheme, meanwhile, the hybrid HWENO scheme is more compact for only immediate neighbor information is needed in the spatial reconstruction.}
\begin{table}
\begin{center}
\caption{2D-Burgers' equation: initial data
$u(x,y,0)=0.5+sin(\pi (x+y)/2)$. Hybrid HWENO and WENO schemes. $T=0.5/\pi$. $L^1$ and $L^\infty$ errors and orders. Uniform meshes with $N_x\times N_y$ cells.}
\medskip
\begin{tabular} {|c|c|c|c|c|c|c|c|c|}
\hline
& \multicolumn{4}{c|}{Hybrid HWENO  scheme} & \multicolumn{4}{c|}{WENO scheme}\\
\hline\hline
  $N_x\times N_y$ cells &$ L^1$ error &  order & $L^\infty$error &  order &$ L^1$ error &  order & $L^\infty$ error &order\\   \hline
$  40\times  40$ &     7.07E-05 &        &     6.32E-04 &       								&     8.20E-05 &        &     6.74E-04 &        \\ \hline
$  80\times  80$ &     3.95E-06 &       4.16 &     4.28E-05 &       3.88 								&     4.06E-06 &       4.34 &     3.91E-05 &       4.11 \\ \hline
$ 120\times 120$ &     7.31E-07 &       4.16 &     7.62E-06 &       4.26 								&     6.30E-07 &       4.60 &     5.67E-06 &       4.76 \\ \hline
$ 160\times 160$ &     2.19E-07 &       4.19 &     2.30E-06 &       4.16 								&     1.66E-07 &       4.64 &     1.42E-06 &       4.81 \\ \hline
$ 200\times 200$ &     8.67E-08 &       4.15 &     8.96E-07 &       4.23 								&     5.65E-08 &       4.82 &     4.96E-07 &       4.72 \\ \hline
$ 240\times 240$ &     4.07E-08 &       4.14 &     4.19E-07 &       4.17 								&     2.28E-08 &       4.99 &     2.22E-07 &       4.40 \\ \hline
\end{tabular}
\label{tburgers2d}
\end{center}
\end{table}
\smallskip

\begin{figure}
 \centerline{
\psfig{file=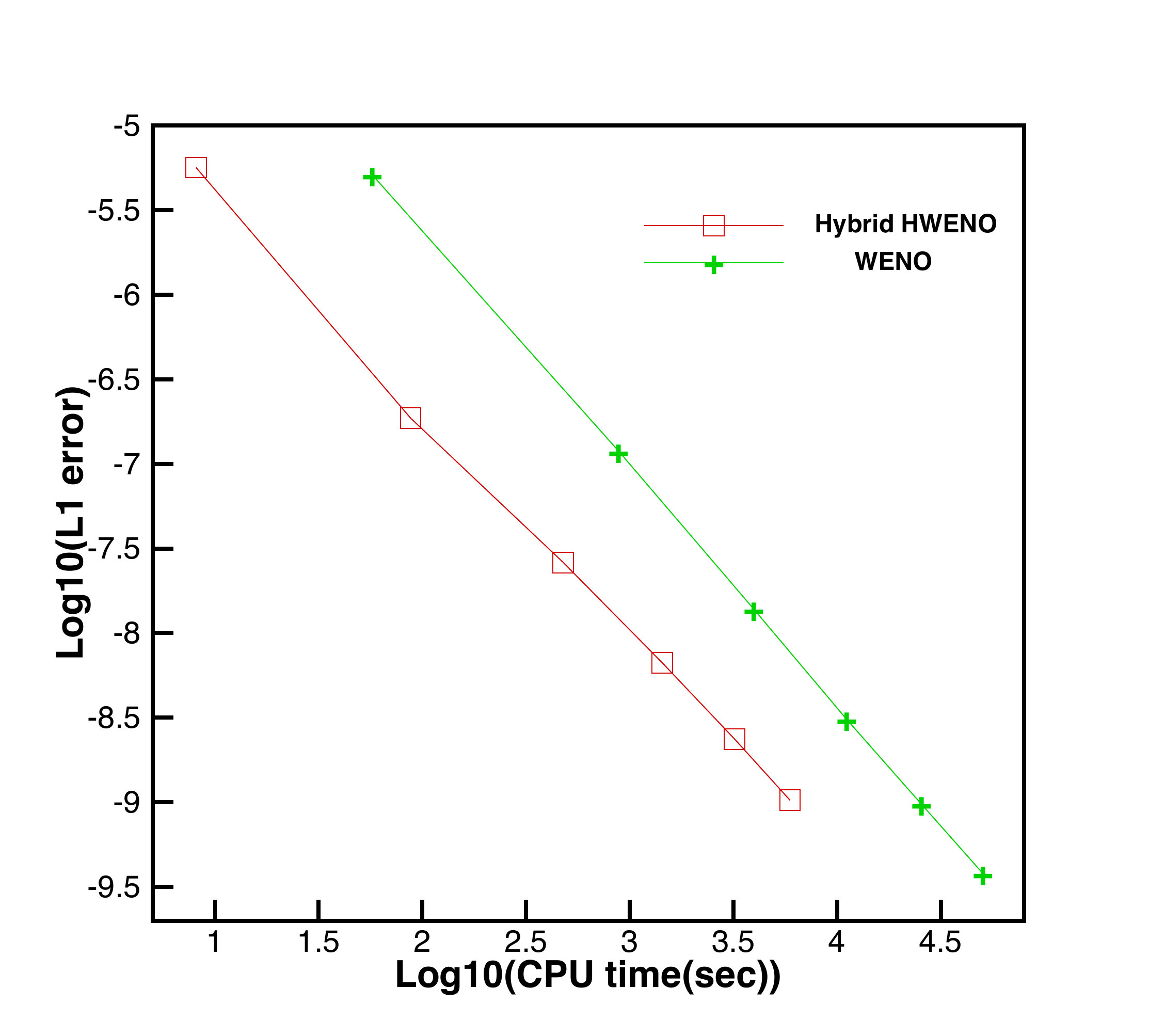,width=2.5 in}
\psfig{file=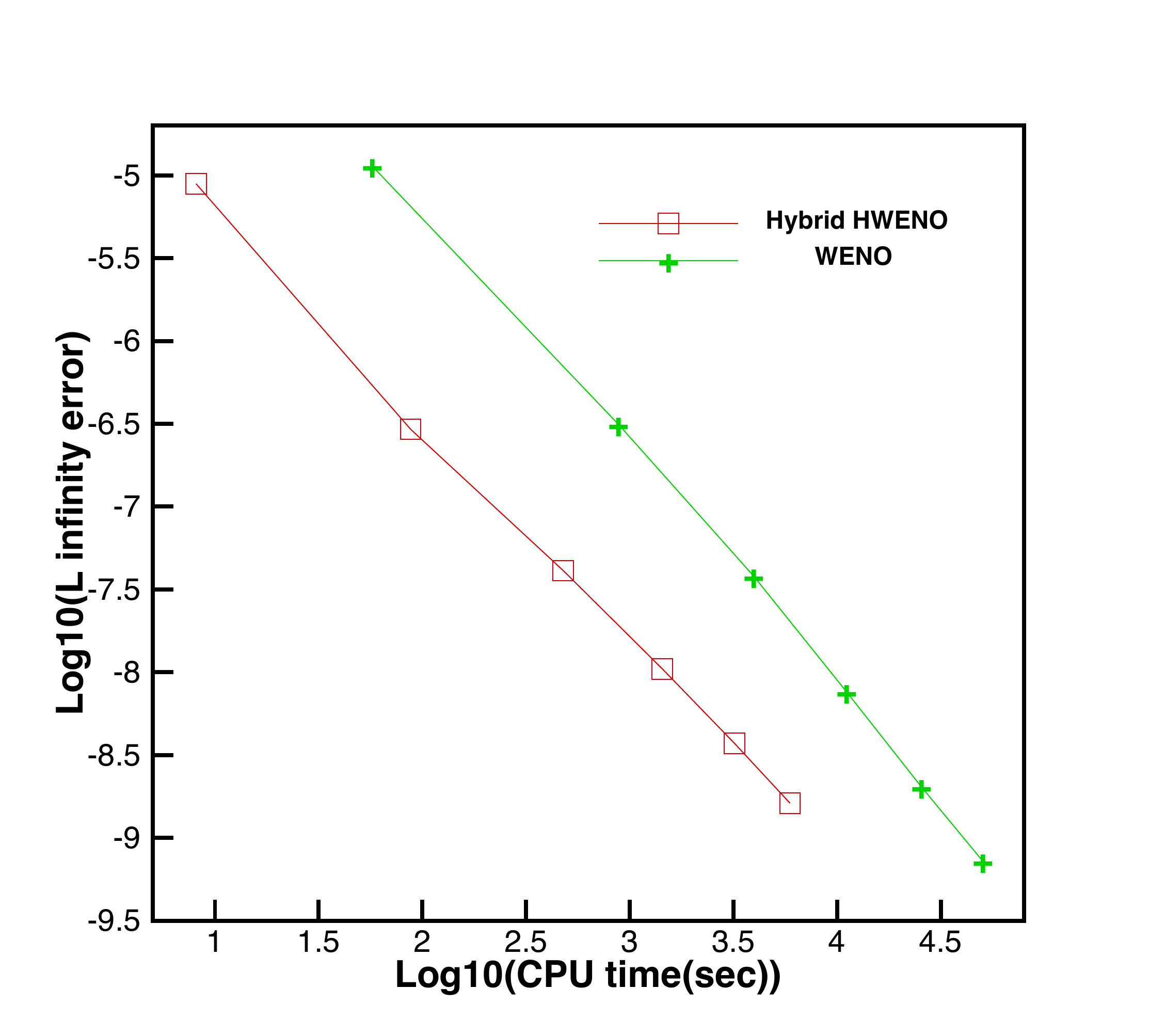,width=2.5 in}}
 \caption{2D-Euler equations: initial data
$\rho(x,y,0)=1+0.2\sin(\pi(x+y))$, $\mu(x,y,0)=1$,
$\nu(x,y,0)=1$ and $p(x,y,0)=1$. $T=2$. Computing times and errors. Rectangle signs and a red solid  line denote the results of Hybrid HWENO scheme;  plus signs and a green solid  line denote the results of WENO scheme.}
\label{FEuler2d_smooth}
\end{figure}

\noindent {\bf Example 3.4.} We consider  two dimensional Euler equations:
\begin{equation}
\label{euler2}
 \frac{\partial}{\partial t}
 \left(
 \begin{array}{c}
 \rho \\
 \rho \mu \\
 \rho \nu \\
 E
 \end{array} \right )
 +
 \frac{\partial}{\partial x}
 \left (
\begin{array}{c}
\rho \mu \\
\rho \mu^{2}+p \\
\rho \mu \nu \\
\mu(E+p)
\end{array}
 \right )
 +
 \frac{\partial}{\partial y}
 \left (
\begin{array}{c}
\rho \nu \\
\rho \mu \nu \\
\rho \nu^{2}+p \\
\nu(E+p)
\end{array}
 \right )=0.
\end{equation}
Here  $\rho$  is the density, $(\mu,\nu)$ is the velocity, $E$ is the total energy, and $p$ is the pressure. The initial conditions are taken as $\rho(x,y,0)=1+0.2\sin(\pi(x+y))$, $\mu(x,y,0)=1$, $\nu(x,y,0)=1$, $p(x,y,0)=1$ and $\gamma=1.4$. The computing domain is $(x,y)\in [0,2] \times [0, 2]$. Periodic boundary conditions are applied in $x$ and $y$ directions. The exact solution of $\rho$ is $\rho(x,y,t)=1+0.2\sin(\pi(x+y-2t))$, $\mu(x,y,0)=1$, $\nu(x,y,0)=1$, $p(x,y,0)=1$ and the computing time is up to $T=2$.
\Red{Table \ref{tEluer2d} gives the numerical errors and  orders of the density for the hybrid HWENO and WENO schemes, and we can know  both two schemes achieve the designed fourth and fifth order accuracy, respectively. In addition, we also present their numerical errors against CPU times in Figure \ref{FEuler2d_smooth}, which shows  Hybrid HWENO scheme has  higher efficiency than WENO scheme}
\begin{table}
\begin{center}
\caption{2D-Euler equations: initial data $\rho(x,y,0)=1+0.2\sin(\pi(x+y))$, $\mu(x,y,0)=1$,
$\nu(x,y,0)=1$ and $p(x,y,0)=1$. Hybrid HWENO and WENO schemes. $T=2$. $L^1$ and $L^\infty$ errors and orders. Uniform meshes with $N_x\times N_y$ cells.}
\medskip
\begin{tabular} {|c|c|c|c|c|c|c|c|c|}
\hline
& \multicolumn{4}{c|}{Hybrid HWENO  scheme} & \multicolumn{4}{c|}{WENO scheme}\\
\hline\hline
$N_x\times N_y$ cells &$ L^1$ error &  order & $L^\infty$error &  order &$ L^1$ error &  order & $L^\infty$ error &order\\   \hline
$  40\times  40$ &     5.66E-06 &            &     8.90E-06 &          								&     5.11E-06 &           &     1.14E-05 &           \\ \hline
$  80\times  80$ &     1.86E-07 &       4.92 &     2.93E-07 &       4.93								&     1.19E-07 &       5.43 &     3.12E-07 &       5.18 \\ \hline
$ 120\times 120$ &     2.61E-08 &       4.85 &     4.10E-08 &       4.85								&     1.38E-08 &       5.31 &     3.79E-08 &       5.20 \\ \hline
$ 160\times 160$ &     6.66E-09 &       4.75 &     1.05E-08 &       4.75								&     3.09E-09 &       5.21 &     7.58E-09 &       5.59 \\ \hline
$ 200\times 200$ &     2.36E-09 &       4.65 &     3.71E-09 &       4.65								&     9.73E-10 &       5.18 &     2.03E-09 &       5.91 \\ \hline
$ 240\times 240$ &     1.03E-09 &       4.56 &     1.62E-09 &       4.55								&     3.79E-10 &       5.17 &     7.21E-10 &       5.67 \\ \hline
\end{tabular}
\label{tEluer2d}
\end{center}
\end{table}
\smallskip

\begin{figure}
 \centerline{
\psfig{file=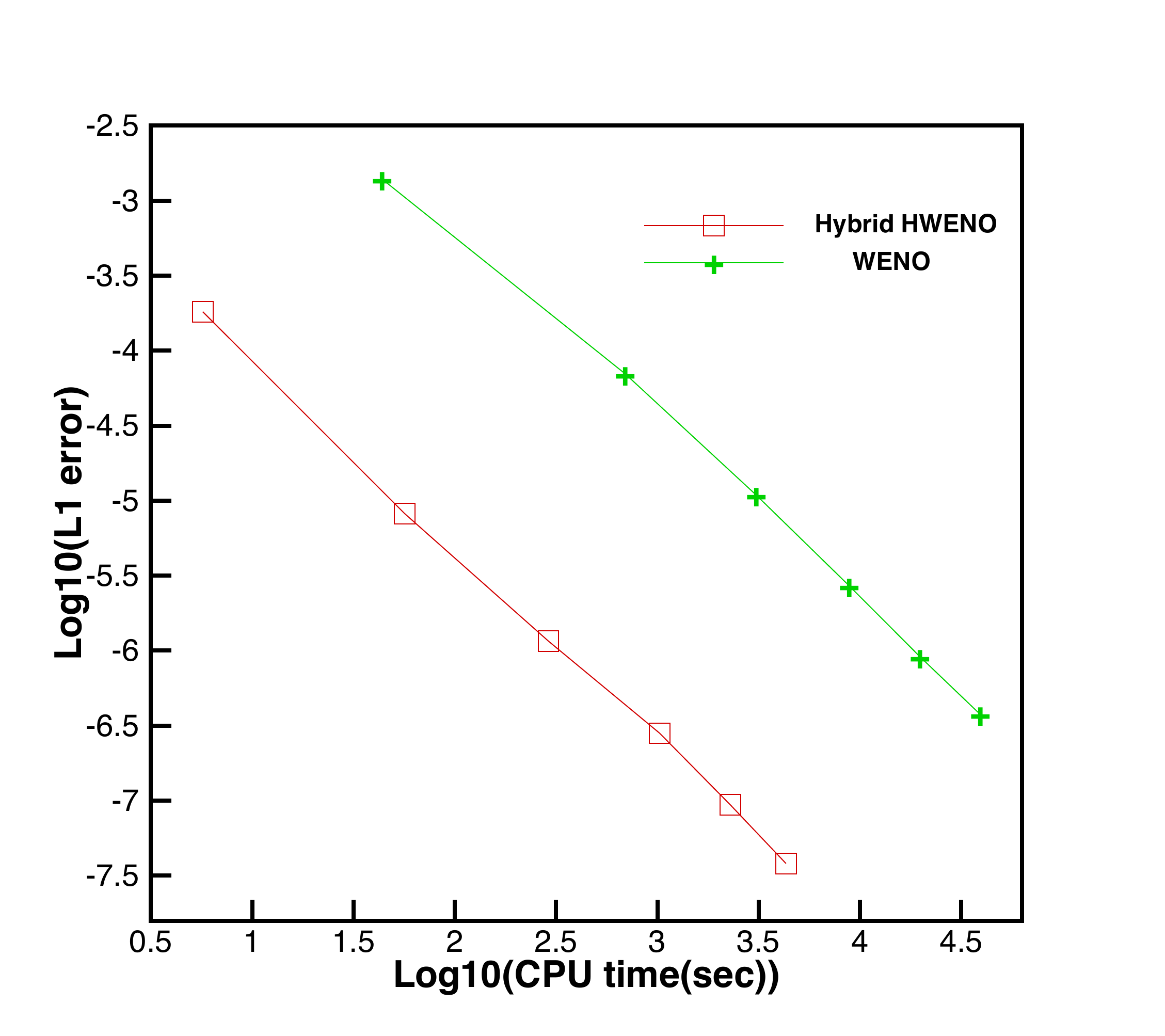,width=2.5 in}
\psfig{file=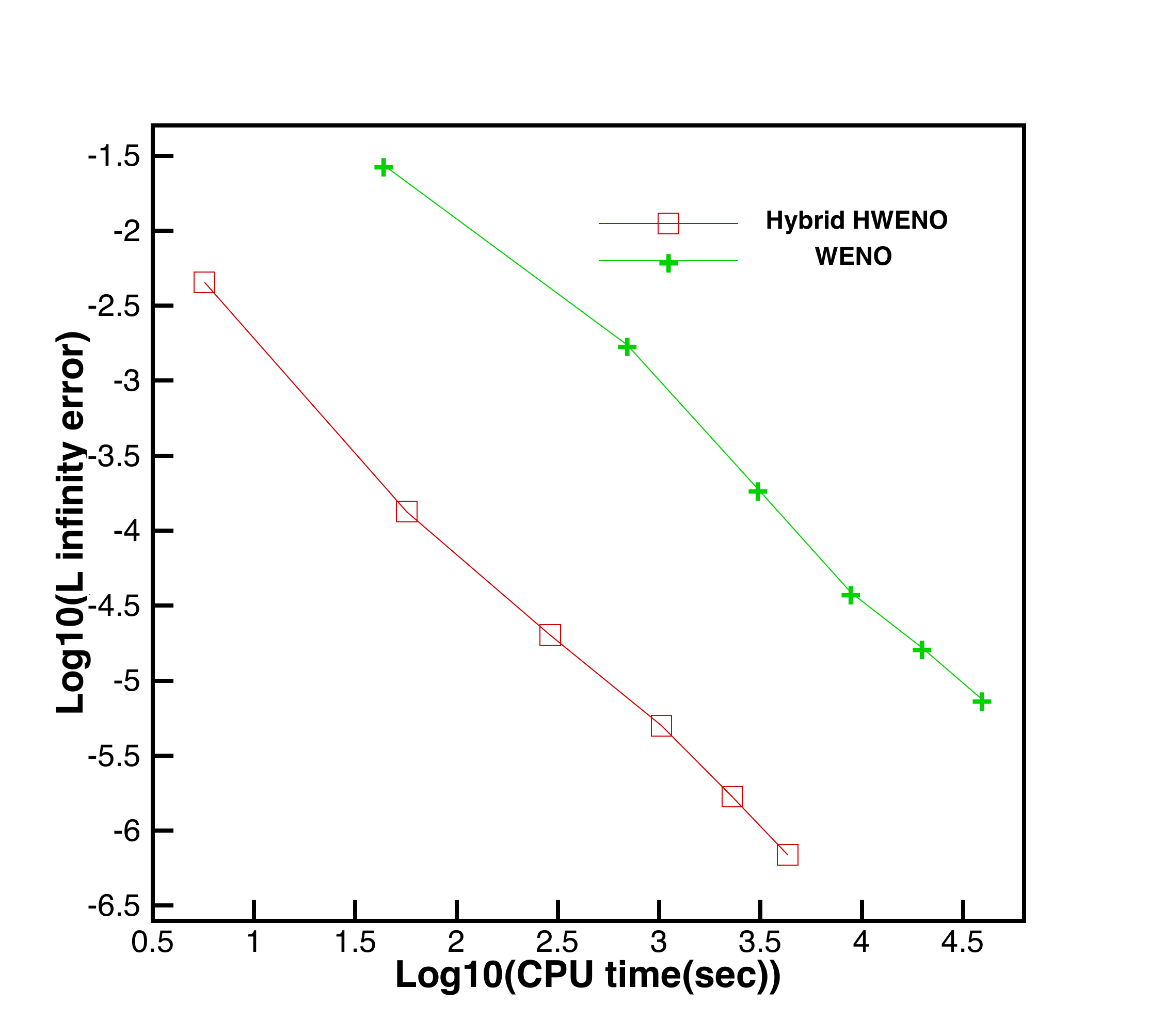,width=2.5 in}}
 \caption{Isentropic vortex test. $T=10$. Computing times and errors. Rectangle signs and a red solid  line denote the results of Hybrid HWENO scheme; plus signs and a green solid  line denote the results of WENO scheme.}
\label{FEuler2d_vortex}
\end{figure}

\Red{\noindent {\bf Example 3.5.} Isentropic vortex test \cite{s2} for
two dimensional Euler equations. A isentropic vortex is added to the mean flow $(\rho, \mu, \nu, p,\gamma)$ $=$ $(1,1,1,1,1.4)$ with perturbations in ($\mu, \nu$), the temperature $T=p/\rho$ and no
perturbation in the entropy $S=p/\rho^\gamma$, so that the initial
conditions are
\begin{equation*}
\begin{split}
\rho&=\left[1-\frac{(\gamma-1)\epsilon^2}{8\gamma\pi^2}e^{1-r^2}\right]^\frac{1}{\gamma-1}, \quad p=\rho^\gamma\\
\mu&=1-\frac{\epsilon y}{2\pi}e^{\frac{1-r^2}{2}},\quad\nu=1+\frac{\epsilon x}{2\pi}e^{\frac{1-r^2}{2},}
\end{split}
\end{equation*}
in which $\epsilon$ represents the vortex strength taken as 5 here and $r^2=x^2+y^2$. The computational domain is $[-5,5]\times[-5,5]$, and periodic boundary conditions are applied in $x$ and $y$ directions. The vortex movement is
aligned with the diagonal of the computational domain, and it returns to the initial position with time period 10. we present the numerical errors and  orders of the density for the hybrid HWENO and WENO schemes in Table \ref{tEluer2dvortex}, then we can see  both two schemes achieve the designed accuracy, respectively. In addition, we also give their numerical errors against CPU times in Figure \ref{FEuler2d_vortex}, which illustrates  Hybrid HWENO scheme has  higher efficiency than WENO scheme with smaller numerical errors and less CPU times.}
\begin{table}
\begin{center}
\caption{Isentropic vortex test. Hybrid HWENO  and WENO schemes. $T=10$. $L^1$ and $L^\infty$ errors and orders. Uniform meshes with $N_x\times N_y$ cells.}
\medskip
\begin{tabular} {|c|c|c|c|c|c|c|c|c|}
\hline
& \multicolumn{4}{c|}{Hybrid HWENO  scheme} & \multicolumn{4}{c|}{WENO scheme}\\
\hline\hline
$N_x\times N_y$ cells &$ L^1$ error &  order & $L^\infty$error &  order &$ L^1$ error &  order & $L^\infty$ error &order\\   \hline
$ 40\times40$ &     1.82E-04 &        &     4.53E-03 &        								 &     1.39E-03 &        &     2.73E-02 &        \\ \hline
$ 80\times80$ &     8.18E-06 &       4.47 &     1.34E-04 &       5.08 								 &     7.00E-05 &       4.32 &     1.73E-03 &       3.98 \\ \hline
$120\times120$&     1.16E-06 &       4.82 &     2.01E-05 &       4.67 								 &     1.09E-05 &       4.58 &     1.89E-04 &       5.46 \\ \hline
$160\times160$&     2.82E-07 &       4.92 &     5.01E-06 &       4.83 								 &     2.70E-06 &       4.85 &     3.86E-05 &       5.52 \\ \hline
$200\times200$&     9.37E-08 &       4.93 &     1.68E-06 &       4.89 								 &     9.12E-07 &       4.87 &     1.66E-05 &       3.78 \\ \hline
$240\times240$&     3.81E-08 &       4.93 &     6.88E-07 &       4.90 								 &     3.77E-07 &       4.85 &     7.51E-06 &       4.36 \\ \hline
\end{tabular}
\label{tEluer2dvortex}
\end{center}
\end{table}
\smallskip

\subsection{Non-smooth tests}

 We present the results of the  hybrid HWENO scheme for the non-smooth tests. \Red{To make comparisons with WENO scheme, we also give the numerical results for some of  one dimensional non-smooth tests.} In addition, we have computed these tests for HWENO scheme likewise, but the results will not be presented here for saving space. Actually, they have similar performance for the problems with discontinuities, but HWENO scheme uses much more computing time. Moreover, we also test the non-smooth problems by the method that we don't modify the first order moments of any cells and  use HWENO reconstruction at the interface  points of each cell for the spatial discretization. Unfortunately, all non-smooth tests have obvious oscillations near discontinuities, even for the one-dimensional Burgers' equation with smaller CFL number, and some tests don't work in this case even though using smaller CFL number, such as Shu-Osher  and two blast waves problems, and so on, which illustrate that the modification for these first order moments in the troubled cells for the hybrid HWENO scheme  is significant to avoid oscillations and keep the hybrid HWENO scheme be robust. In the implementation, there are only a small part of cells in which we need to modify their first order moments for the hybrid HWENO scheme.

\begin{figure}
 \centerline{ \psfig{file=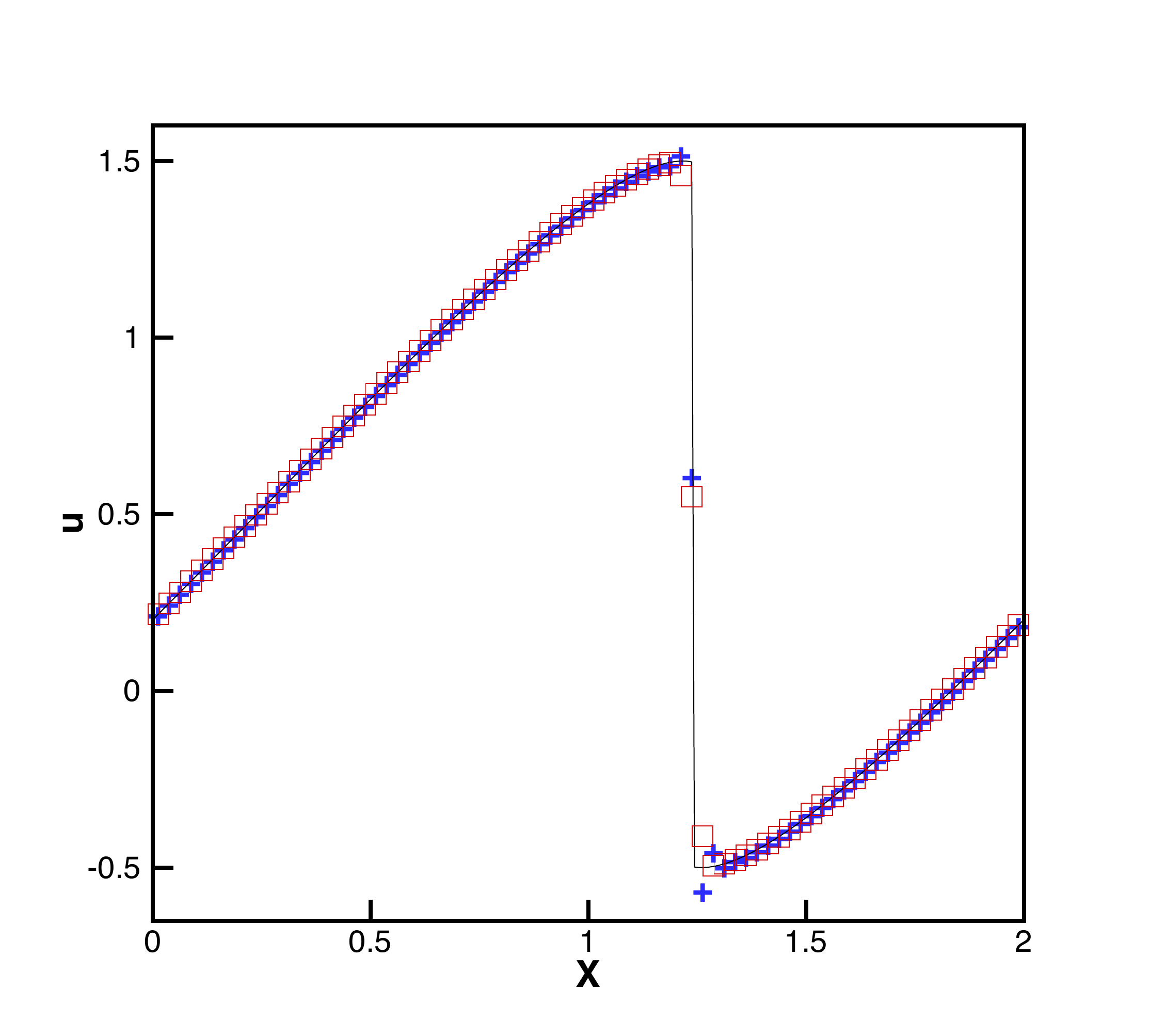,width=3 in}}
 \caption{1D-Burgers' equation: initial data
$u(x,0)=0.5+sin(\pi x)$. $T=1.5/\pi$. Solid line: exact solution; blue plus signs: the results computed by the method that we don't modify the first order moments of any cells and use HWENO reconstruction at the interface  points of each cell for the spatial discretization; red squares: the results of hybrid HWENO scheme. Uniform meshes with 80 cells.}
\label{Fburges1d}
\end{figure}

\noindent{\bf Example 3.6.} We solve the one-dimensional Burgers' equation (\ref{1dbugers}) as in Example 3.1. The same initial and boundary conditions are applied here. The computing time is up to $t=1.5/\pi$, and the solution is discontinuous by this time, then, we plot the numerical solution of the hybrid HWENO scheme and the exact solution in Figure \ref{Fburges1d}, and we also test this example by the method that we don't modify any first order moments and directly use HWENO reconstruction at the interface points of each cell for the spatial discretization. We set the CFL number as 0.45 in this case for the code doesn't work with the original CFL number, and its numerical solution is also presented in  Figure \ref{Fburges1d}. From this figure, we can see that if we don't modify the first order moments in the troubled cell and directly use HWENO procedures at the interface points of each cell for spatial discretization, which would have obvious oscillations even with smaller CFL number. However, we also find that the hybrid HWENO scheme can avoid the non-physical oscillations and has a good resolution, which shows that the modification for the first order moments near the discontinuities is a significant and essential procedure. In addition, there are only 4.52\% cells where we need to modify their first order moments by calculating.
\smallskip

\begin{figure}
 \centerline{\psfig{file=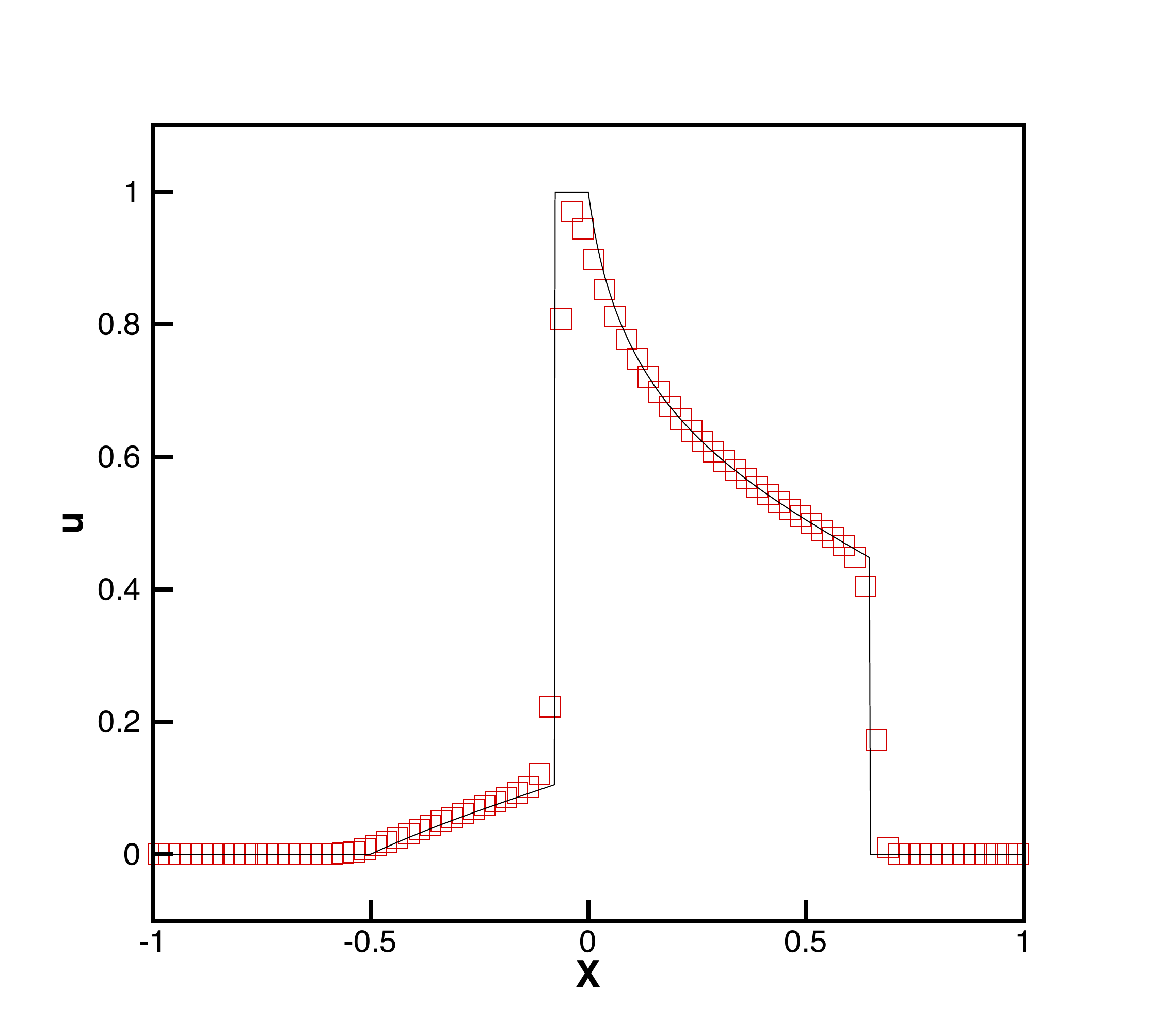,width=3 in}}
 \caption{1D Buckley-Leverett problem: initial data
$u=1$ for $-\frac12\leq x\leq0$ and $u=0$ elsewhere. $T=0.4$. Solid line: exact solution; squares: the results of  the hybrid  HWENO scheme. Uniform meshes with 80 cells.}
\label{FBuk_levert}
\end{figure}

\noindent{\bf Example 3.7.} We now consider a one dimensional nonlinear non-convex scalar Buckley-Leverett problem
\begin{equation}
  u_t+\left( \frac{4u^2}{4u^2+(1-u)^2}\right)_x=0,\quad -1\leq x\leq 1,
\end{equation}
with the initial condition: $u=1$ for $-\frac12\leq x\leq0$ and $u=0$ elsewhere. Inflow and outflow conditions are applied at left and right boundary, respectively, and  the computing time is up to $t=0.4$. The exact solution contains both shock and rarefaction, moreover, some high-order schemes may non-converge to the right entropy solution for this test. We present the numerical results in Figure \ref{FBuk_levert}, and we can find that the hybrid HWENO scheme performs well for capturing the correct entropy solution and has a good resolution.
\smallskip

\begin{figure}
\centering{\psfig{file=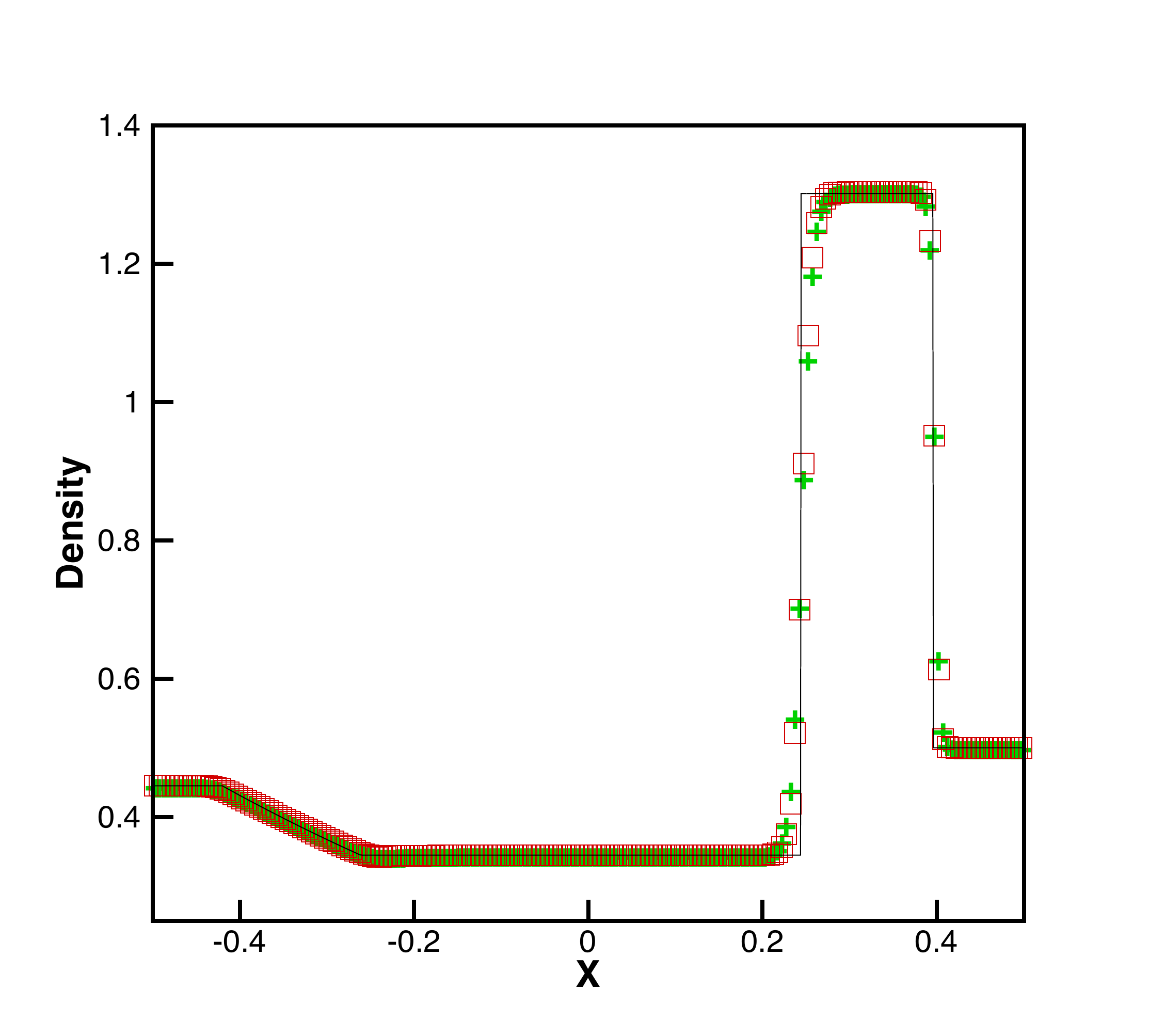,width=2 in}\psfig{file=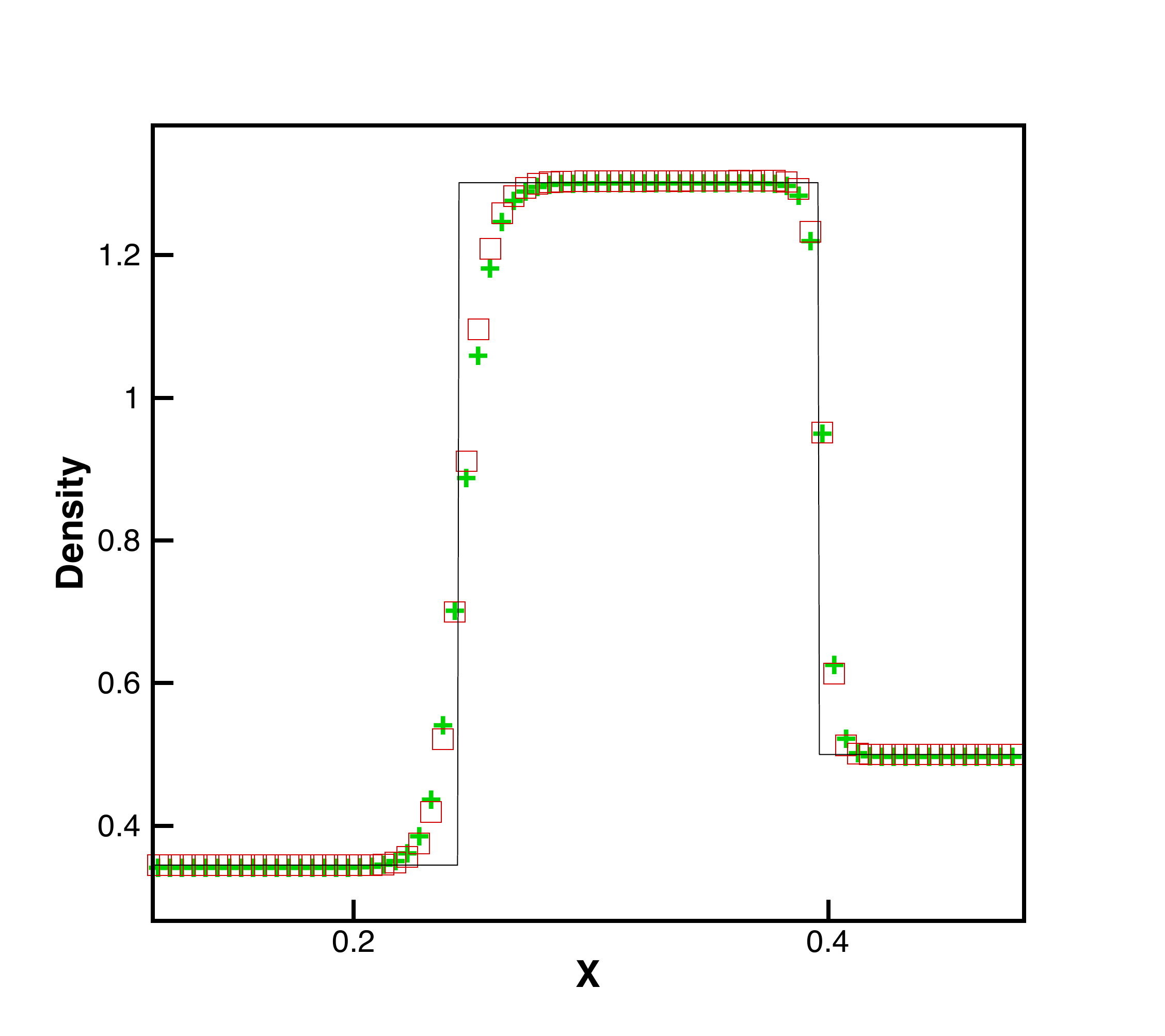,width=2 in}\psfig{file=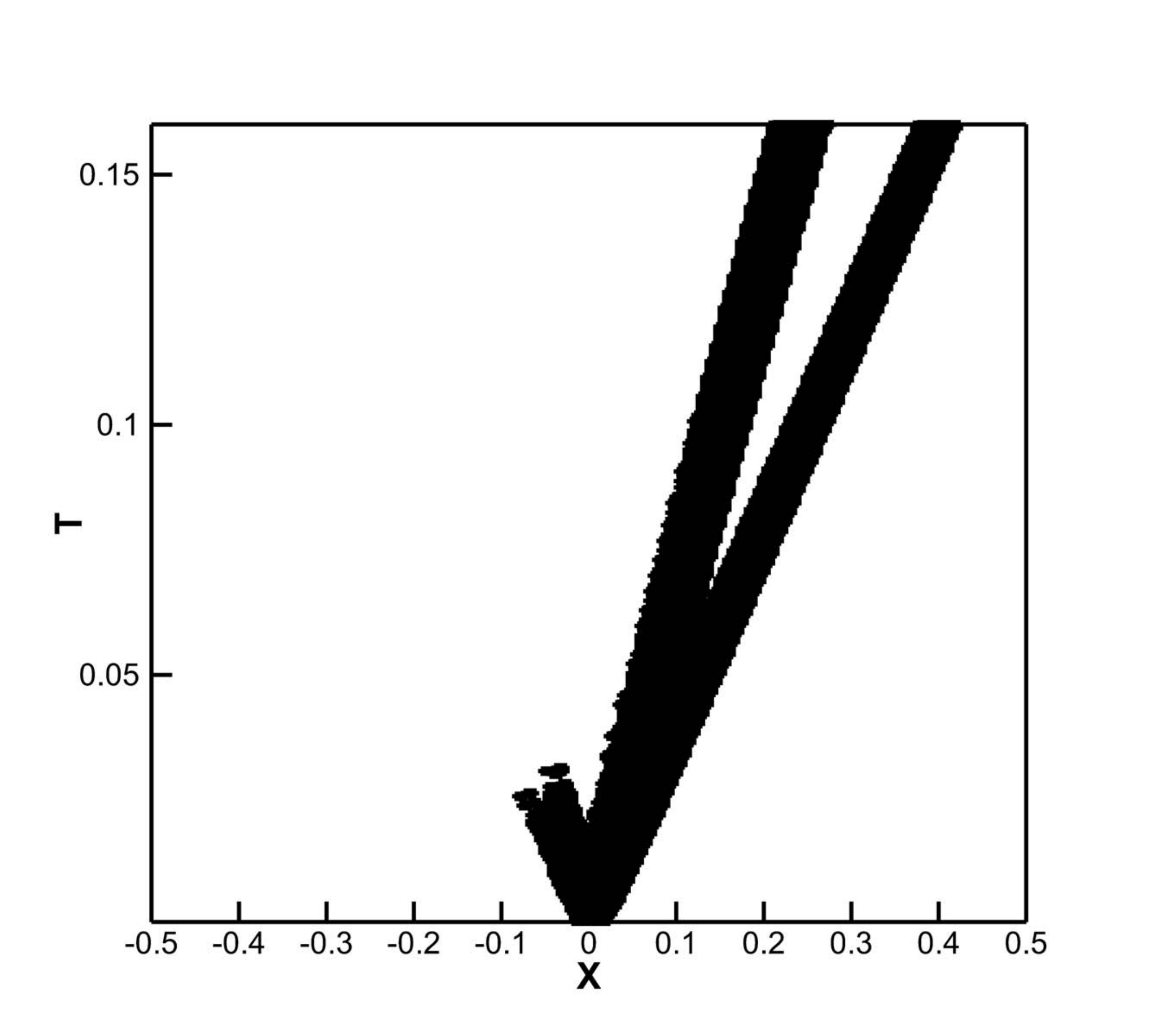,width=2 in}}
 \caption{The Lax problem.  T=0.16. From left to right: density; density zoomed in; the time history of the cells where  we modify the first order moments in the hybrid HWENO scheme. Solid line: the exact solution; red squares: the results of the hybrid  HWENO scheme; green plus signs: the results of WENO scheme. Uniform meshes with 200 cells.}
\label{laxfig}
\end{figure}

\noindent {\bf Example 3.8.} We solve the 1D Euler equations Riemann initial condition for the Lax problem
\begin{equation}
\label{lax} (\rho,\mu,p,\gamma)^T= \left\{
\begin{array}{ll}
(0.445,0.698,3.528,1.4)^T,& x \in [-0.5,0),\\
(0.5,0,0.571,1.4)^T,& x \in [0,0.5].
\end{array}
\right.
\end{equation}
The final computing time is up to $T=0.16$, and we  first present the performances of  the exact solution and the computed density $\rho$ obtained with the hybrid HWENO and WENO schemes  by using 200 uniform cells in Figure \ref{laxfig}. The zoomed in picture and the time history of the cells where  we modify the first order moments in the hybrid HWENO scheme are also given in Figure \ref{laxfig}. In this test case, there are about 10.71\% cells in which we modify their first order moments, which means most regions use linear approximation and have non-decomposition in the characteristic direction, therefore, the hybrid HWENO scheme saves about 62.5\% computational time than the HWENO scheme, meanwhile, the modification in the troubled cells is very important to control  the oscillations. The hybrid HWENO scheme also keeps good resolution, \Red{and the hybrid HWENO and WENO schemes have similar numerical performances, but the hybrid HWENO scheme only needs the immediate neighbor information.}
\smallskip

\begin{figure}
 \centering{\psfig{file=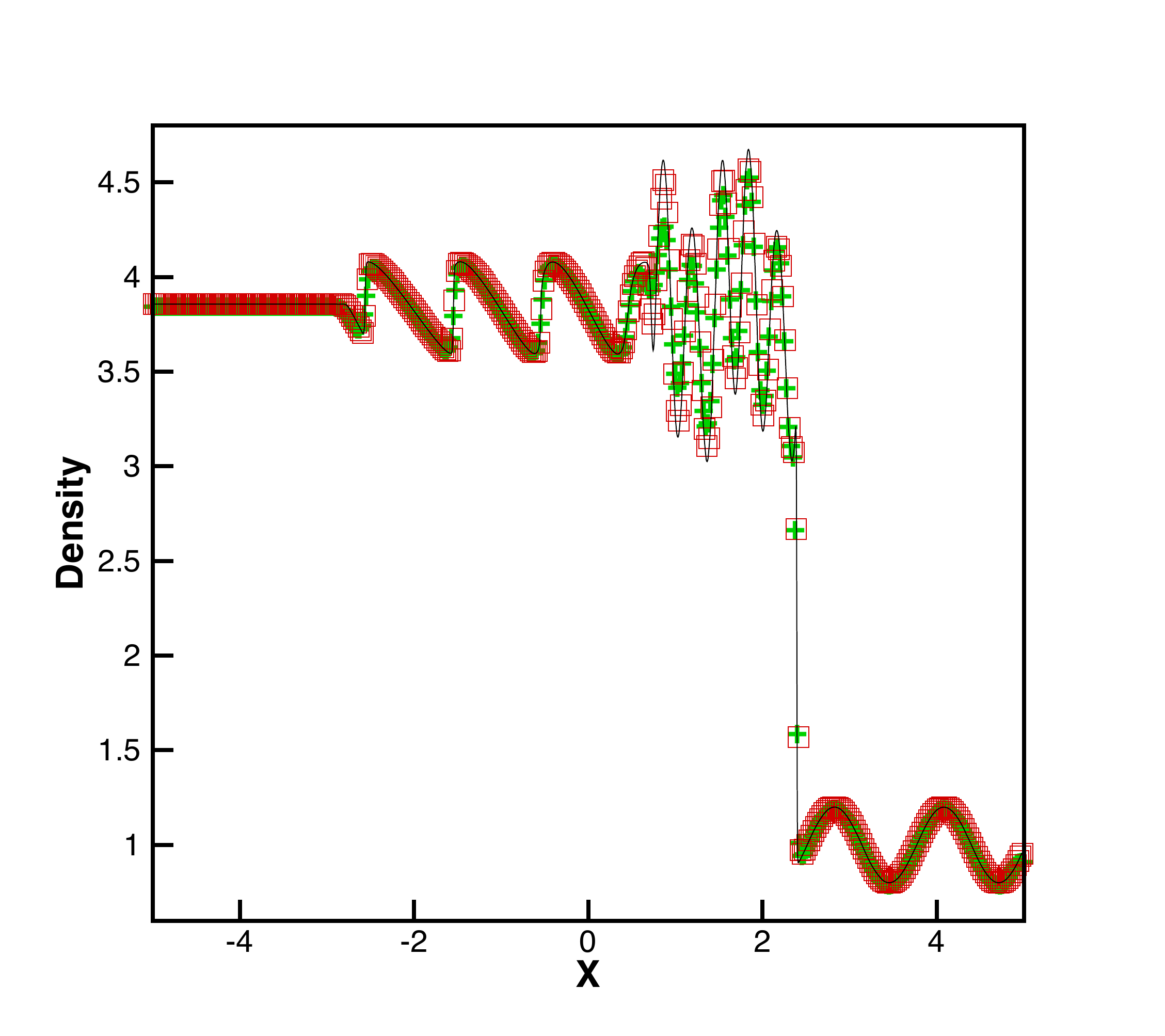,width=2 in}\psfig{file=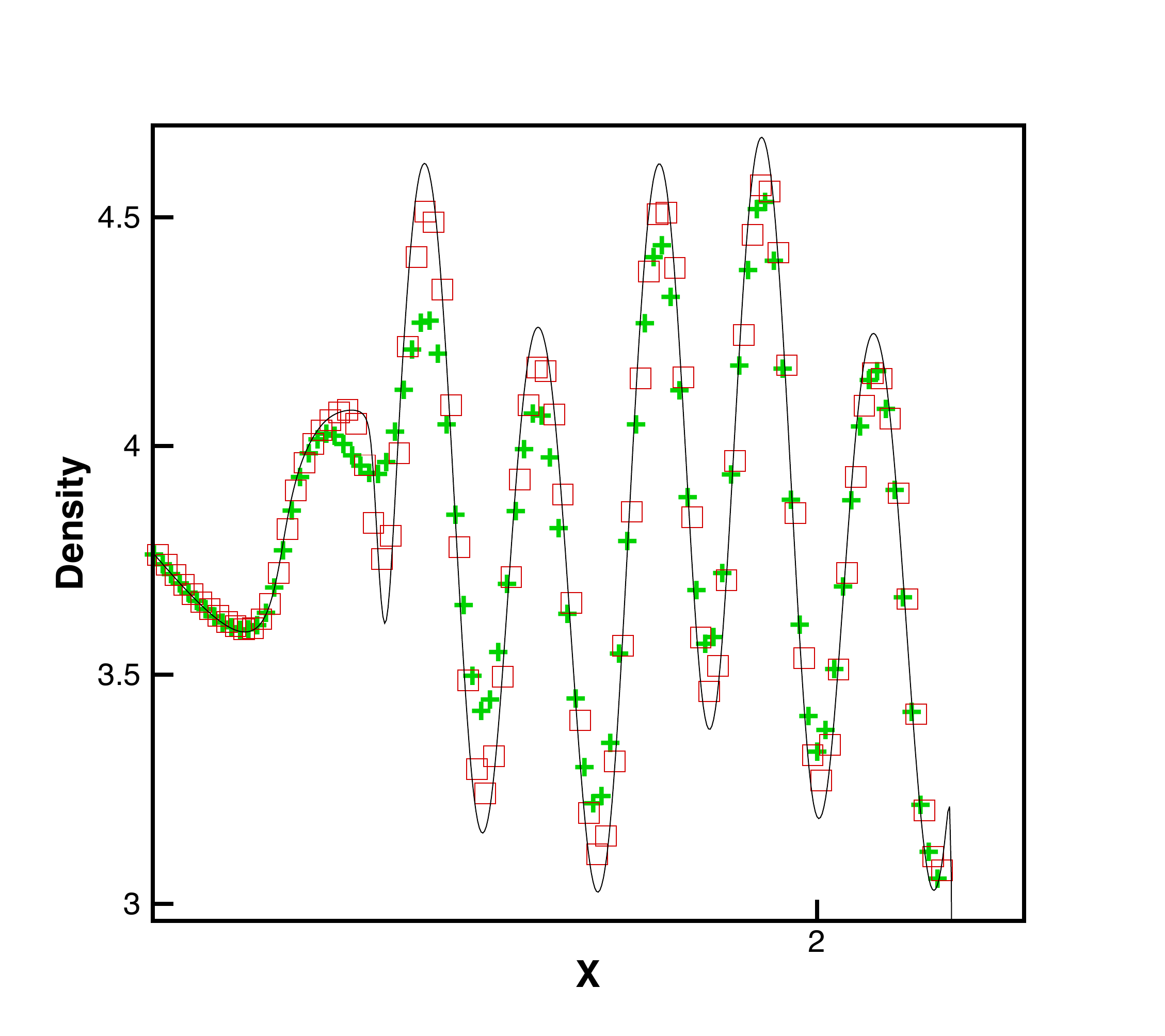,width=2 in}\psfig{file=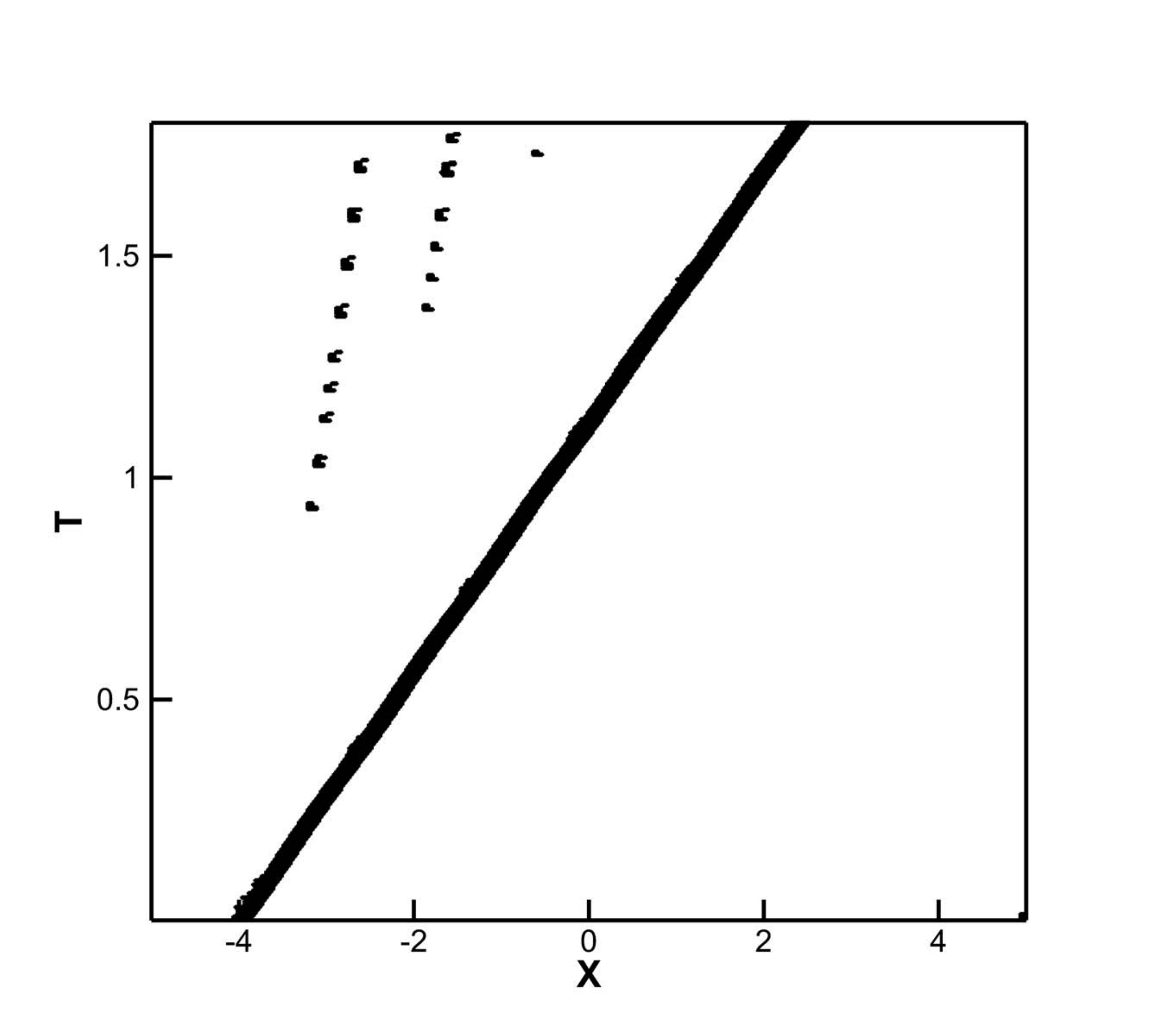,width=2 in}}
\caption{The shock density wave interaction problem. T=1.8. From left to right: density; density zoomed in; the time history of the cells where  we modify the first order moments in the hybrid HWENO scheme. Solid line: the exact solution; red squares: the results of the hybrid  HWENO scheme; green plus signs: the results of WENO scheme. Uniform meshes with 400 cells.}
\label{sin}
\end{figure}

\noindent {\bf Example 3.9.} The Shu-Osher problem, which describes shock interaction with entropy waves \cite{s3}, and the initial condition is
\begin{equation}
\label{ShuOsher} (\rho,\mu,p,\gamma)^T= \left\{
\begin{array}{ll}
(3.857143, 2.629369, 10.333333,1.4)^T,& x \in [-5, -4),\\
(1 + 0.2\sin(5x), 0, 1,1.4)^T,& x \in [-4,5].
\end{array}
\right.
\end{equation}
As we know, when the solutions both contains  shocks and complex smooth region structures, a good scheme would simulate it well. Actually, this test case is a typical example with a moving Mach=3 shock interacting with sine waves in density.
We solve this problem up to $T=1.8$. In Figure \ref{sin}, we present the computed density $\rho$ by the hybrid HWENO scheme and WENO schemes against the referenced "exact" solution, the zoomed in picture and the time history of the cells where  we modify the first order moments in the hybrid HWENO scheme. The referenced "exact" solution is computed by the fifth order finite difference WENO scheme \cite{js} with 2000 grid points. We find that there are only 2.42\% cells where we need to modify their first order moments by calculating, which saves near 66.7\% CPU time than HWENO scheme, but the modification in the troubled cells is very significant to make the scheme be robust. \Red{In addition, we also see that the hybrid HWENO scheme has higher resolution than WENO scheme shown in Figure \ref{sin}. Particularly, the hybrid HWENO scheme only needs three cells while the WENO scheme needs five cells for the spatial reconstruction.}
\smallskip

\begin{figure}
  \centering{\psfig{file=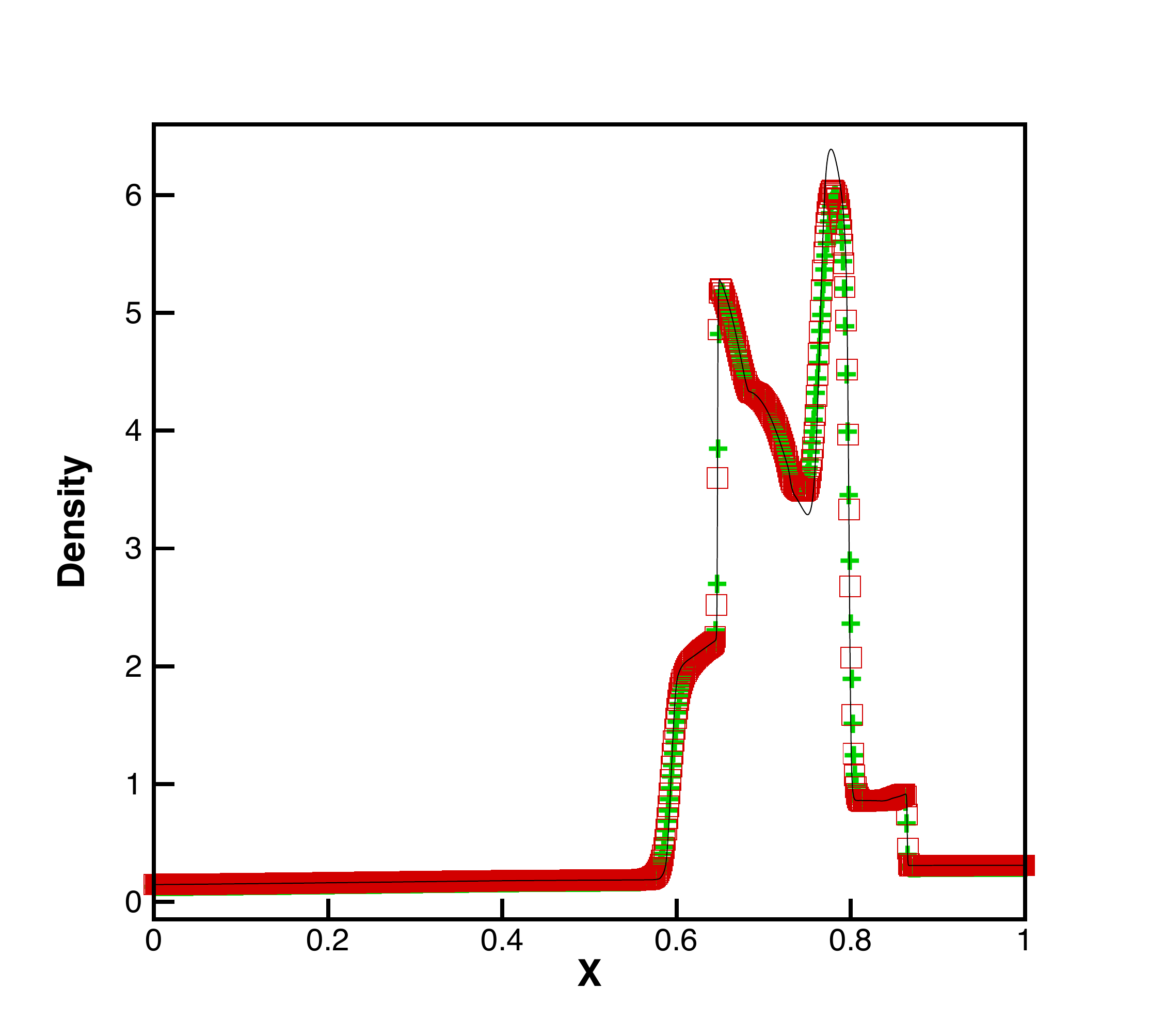,width=2 in}\psfig{file=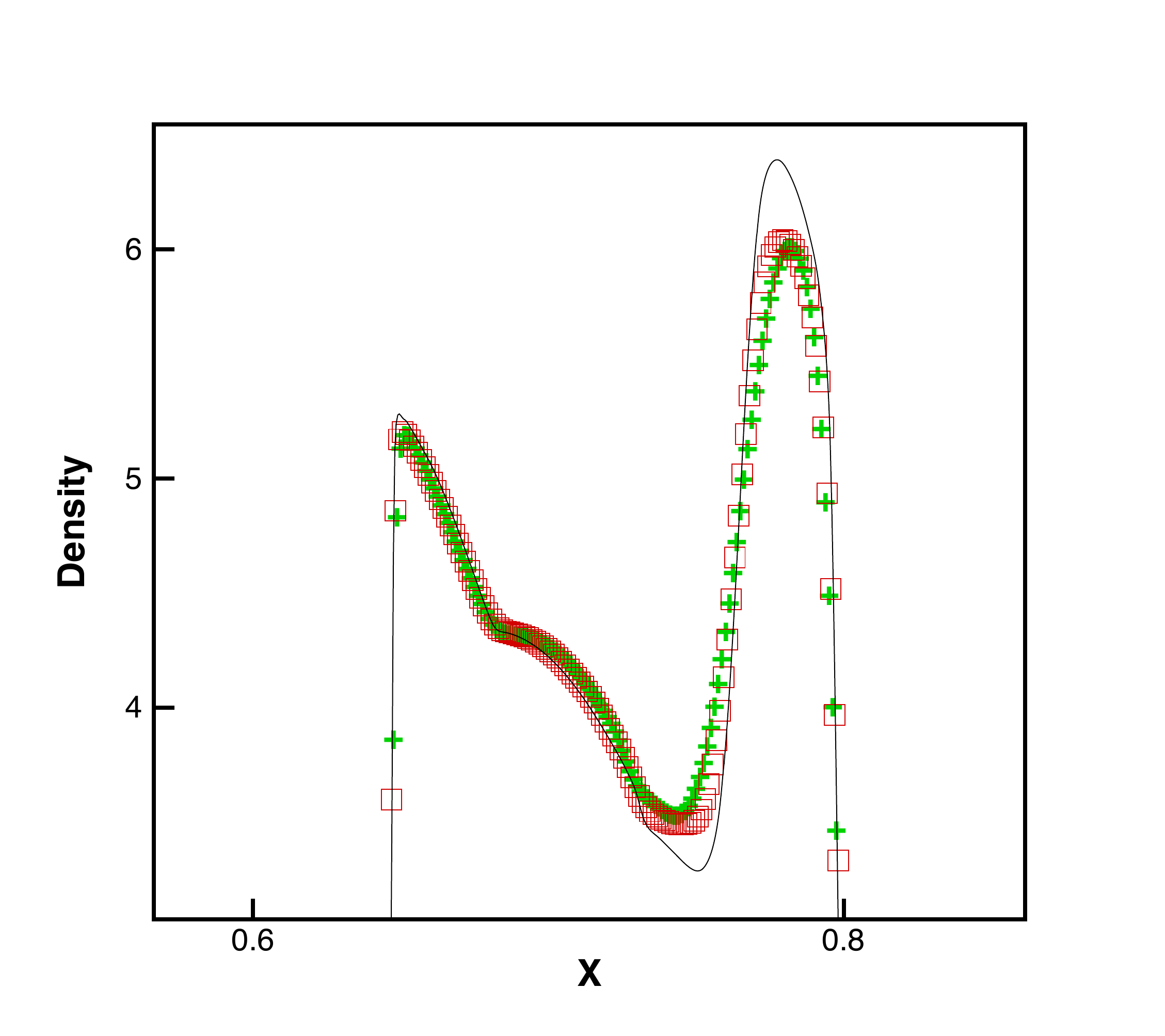,width=2 in}\psfig{file=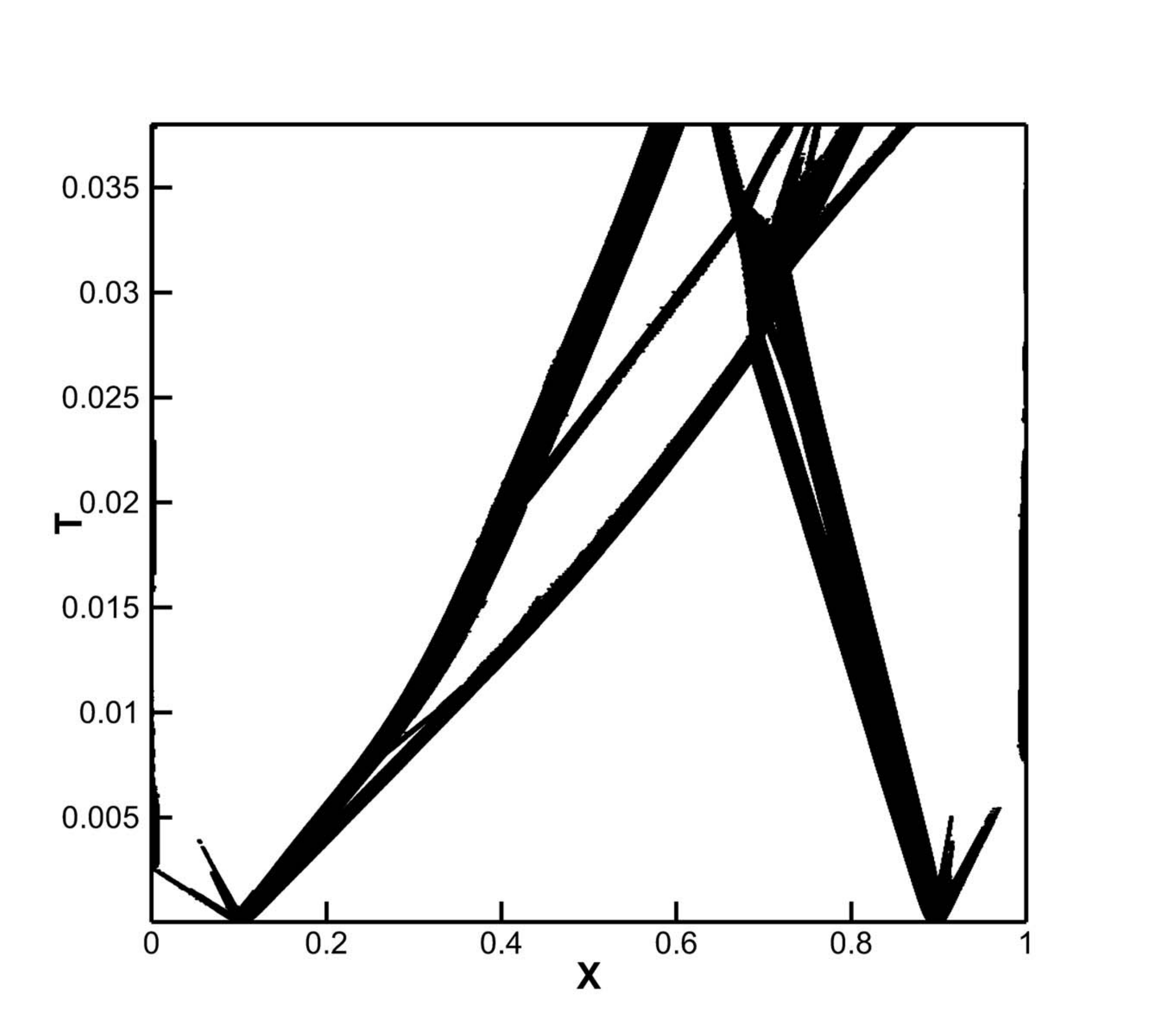,width=2 in}}
\caption{The blast wave problem. T=0.038. From left to right: density; density zoomed in; the time history of the cells where  we modify the first order moments in the hybrid HWENO scheme. Solid line: the exact solution; red squares: the results of the hybrid  HWENO scheme; green plus signs: the results of WENO scheme. Uniform meshes with 800 cells.}
\label{blast}
\end{figure}

\noindent {\bf Example 3.10.} We now consider the interaction of two blast waves, and the initial conditions are:
\begin{equation}
\label{blastwave} (\rho,\mu,p,\gamma)^T= \left\{
\begin{array}{ll}
(1,0,10^3,1.4)^T,& 0<x<0.1,\\
(1,0,10^{-2},1.4)^T,& 0.1<x<0.9,\\
(1,0,10^2,1.4)^T,& 0.9<x<1.
\end{array}
\right.
\end{equation}
The final computing time $T=0.038$ and the reflective condition is applied in two boundaries. In Figure \ref{blast}, we present the computed density by the hybrid HWENO and WENO schemes against the reference "exact" solution, the zoomed in picture  and the time history of the cells where  we modify the first order moments in the hybrid HWENO scheme. Particularly, the reference "exact" solution is a converged solution computed by the fifth order finite difference WENO scheme \cite{js} with 2000 grid points. In the implementation, we find that the hybrid HWENO scheme saves about 58.5\% computational time as there are almost 11.31\% cells where we need to modify their first order moments. The modification for the first order moment in the troubled cells makes the hybrid HWENO scheme be robust, meanwhile, \Red{the hybrid HWENO scheme has higher resolution than WENO scheme.}
\smallskip

\begin{figure}
 \centerline{\psfig{file=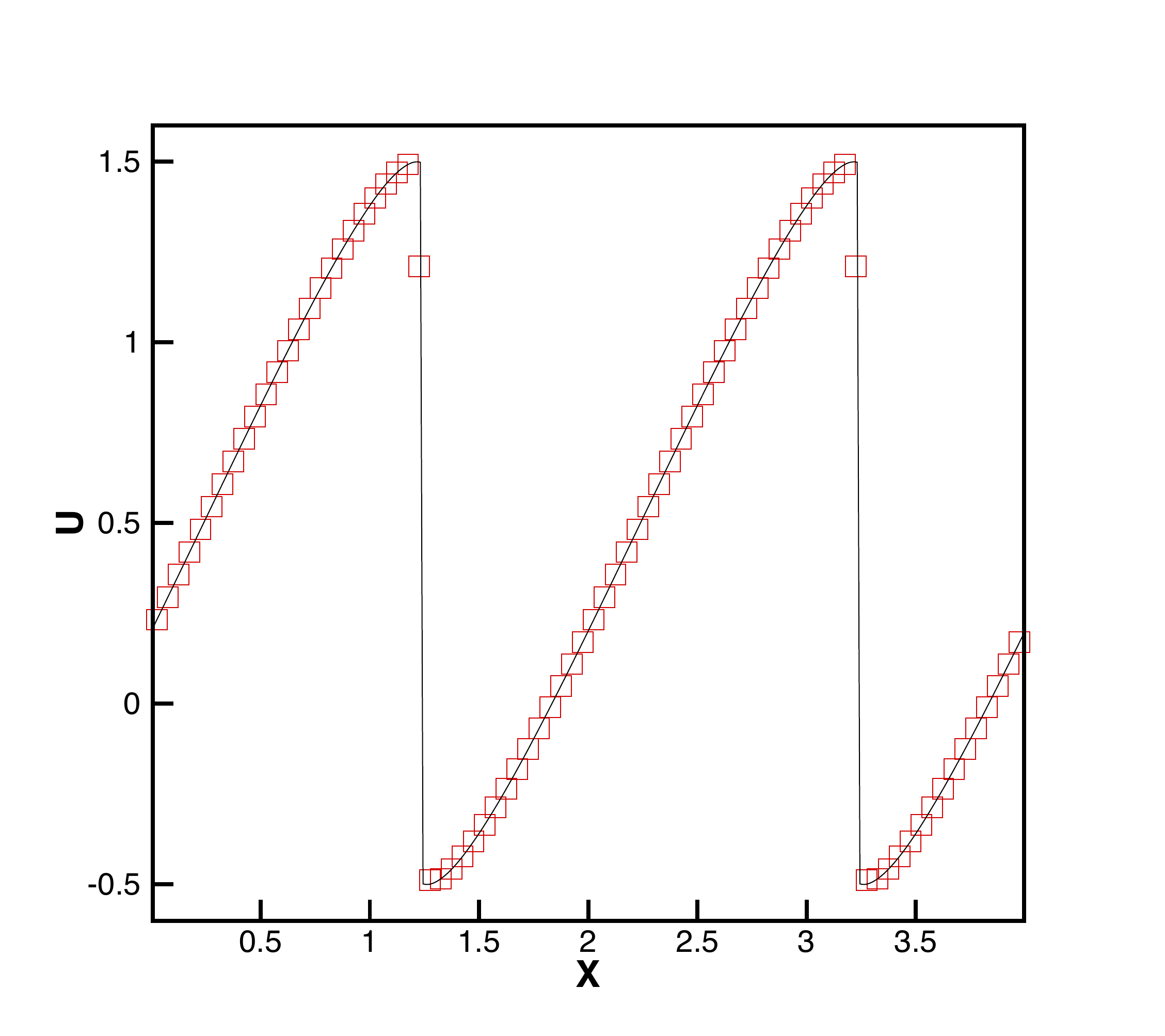,width=2.5 in}\psfig{file=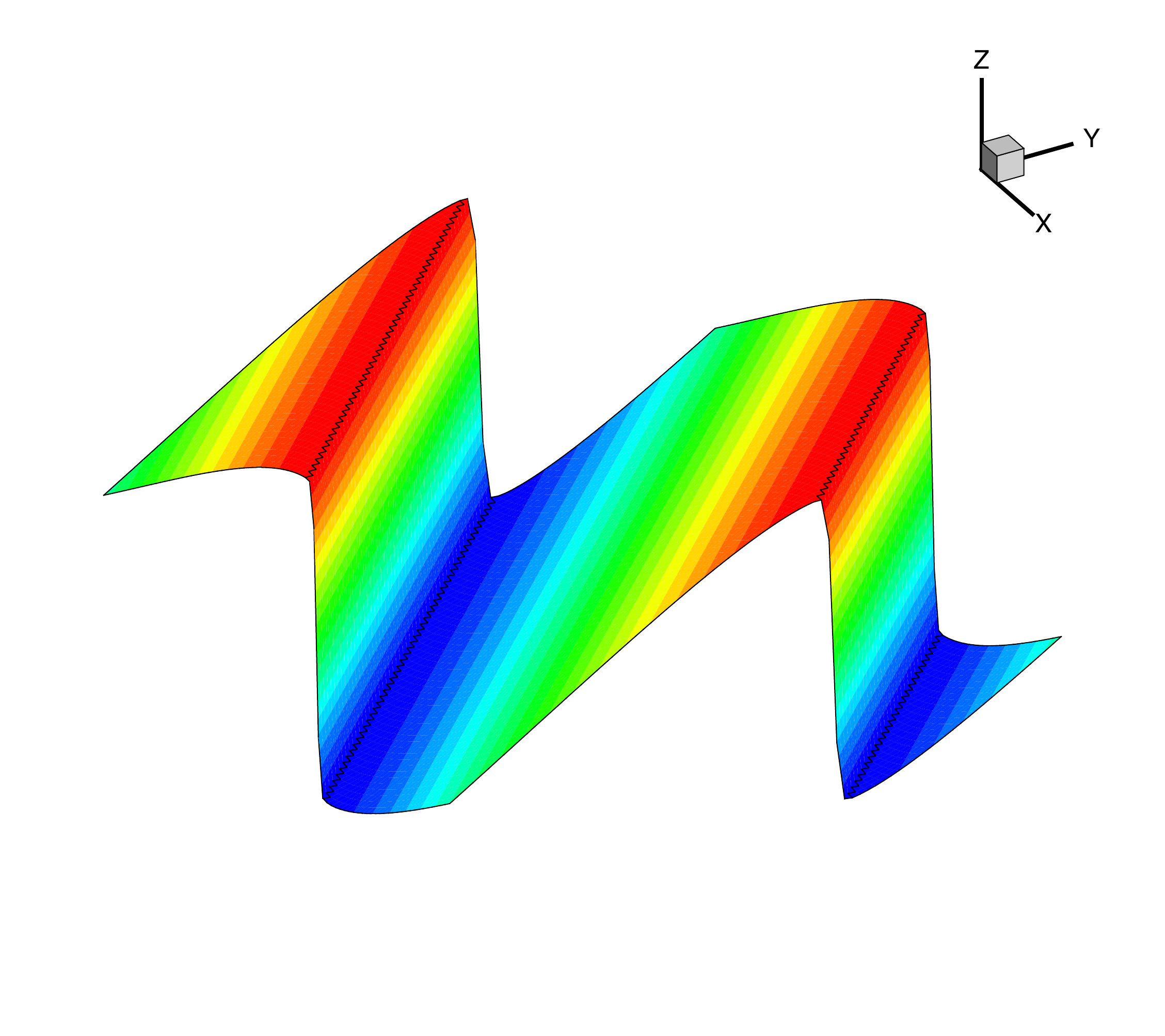,width=2.5 in}}
 \caption{2D-Burgers' equation: initial data
$u(x,y,0)=0.5+sin(\pi (x+y)/2)$. $T=1.5/\pi$. From left to right, the numerical solution at $x=y$ and its surface. Solid line: exact solution; squares: the results of  the hybrid HWENO scheme. Uniform meshes with $80\times80$ cells.  }
\label{Fburges2d}
\end{figure}

\noindent{\bf Example 3.11.} We consider the two-dimensional Burgers' equation (\ref{2dbugers}) seen in Example 3.3 with the same initial and boundary conditions here. The finial computing time is $t=1.5/\pi$, and the solution is discontinuous, then, we plot the numerical solution and the exact solution in Figure \ref{Fburges2d} with $80\times 80$ uniform meshes. Similarly as in one dimensional Burgers' equation, the hybrid HWENO scheme has good resolution nearby discontinuities.
\smallskip

\begin{figure}
\centerline{\psfig{file=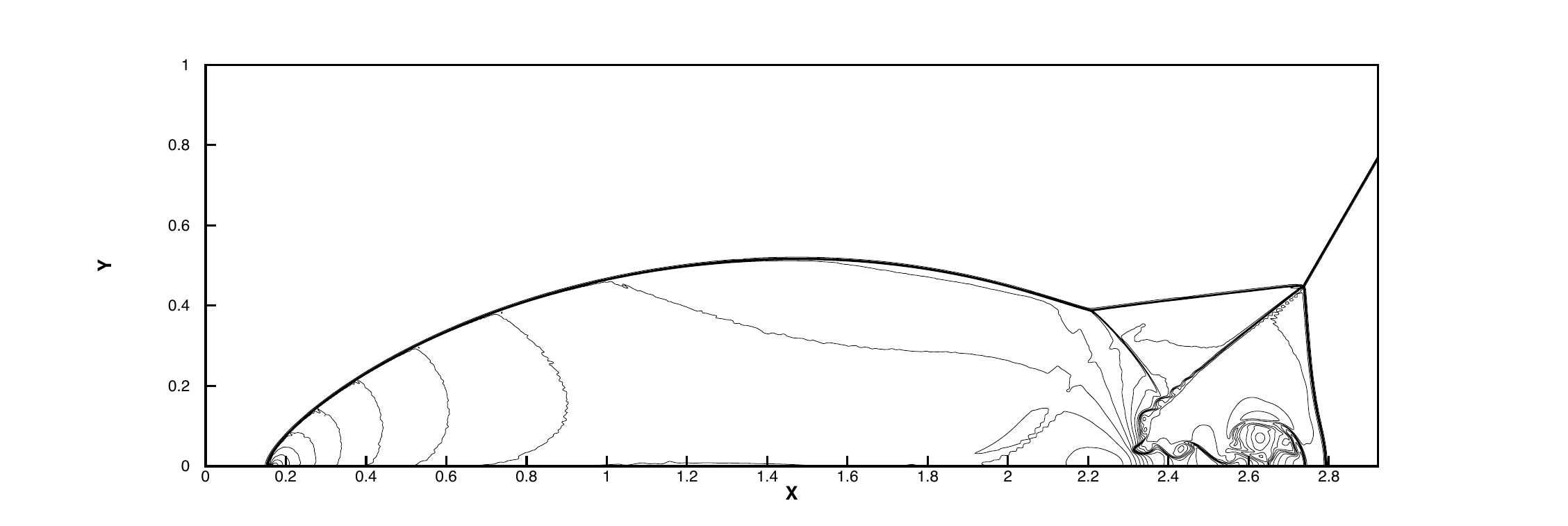,width=6 in}}
\centerline{\psfig{file=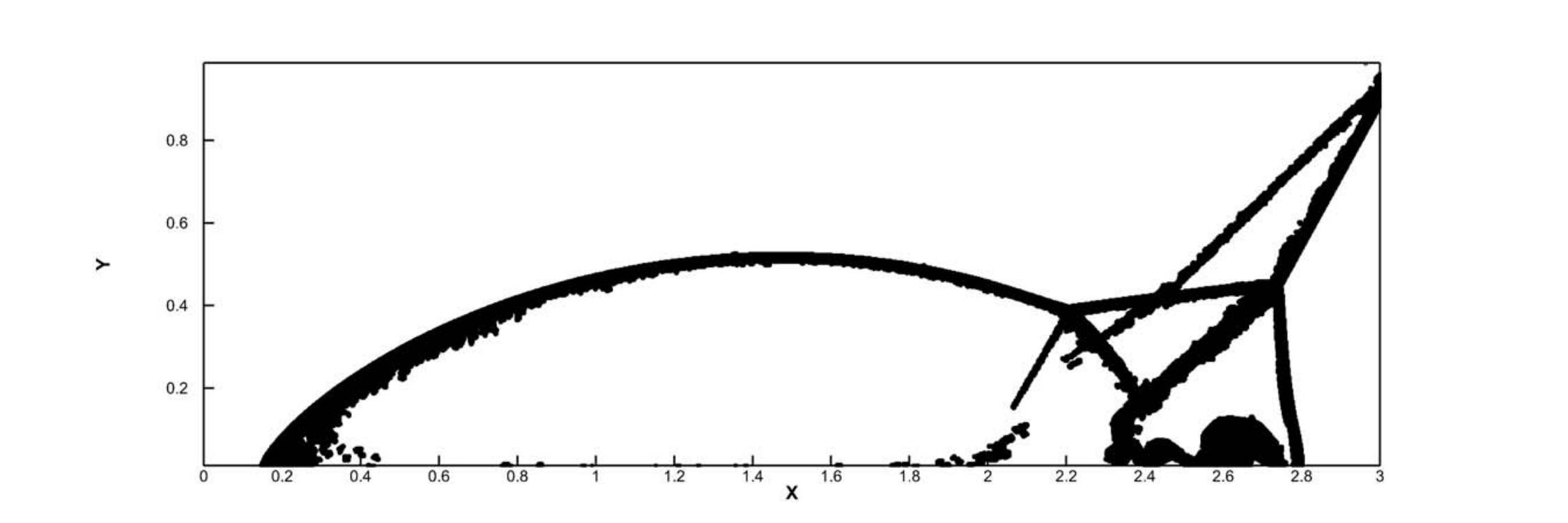,width=6 in}}
\centerline{\psfig{file=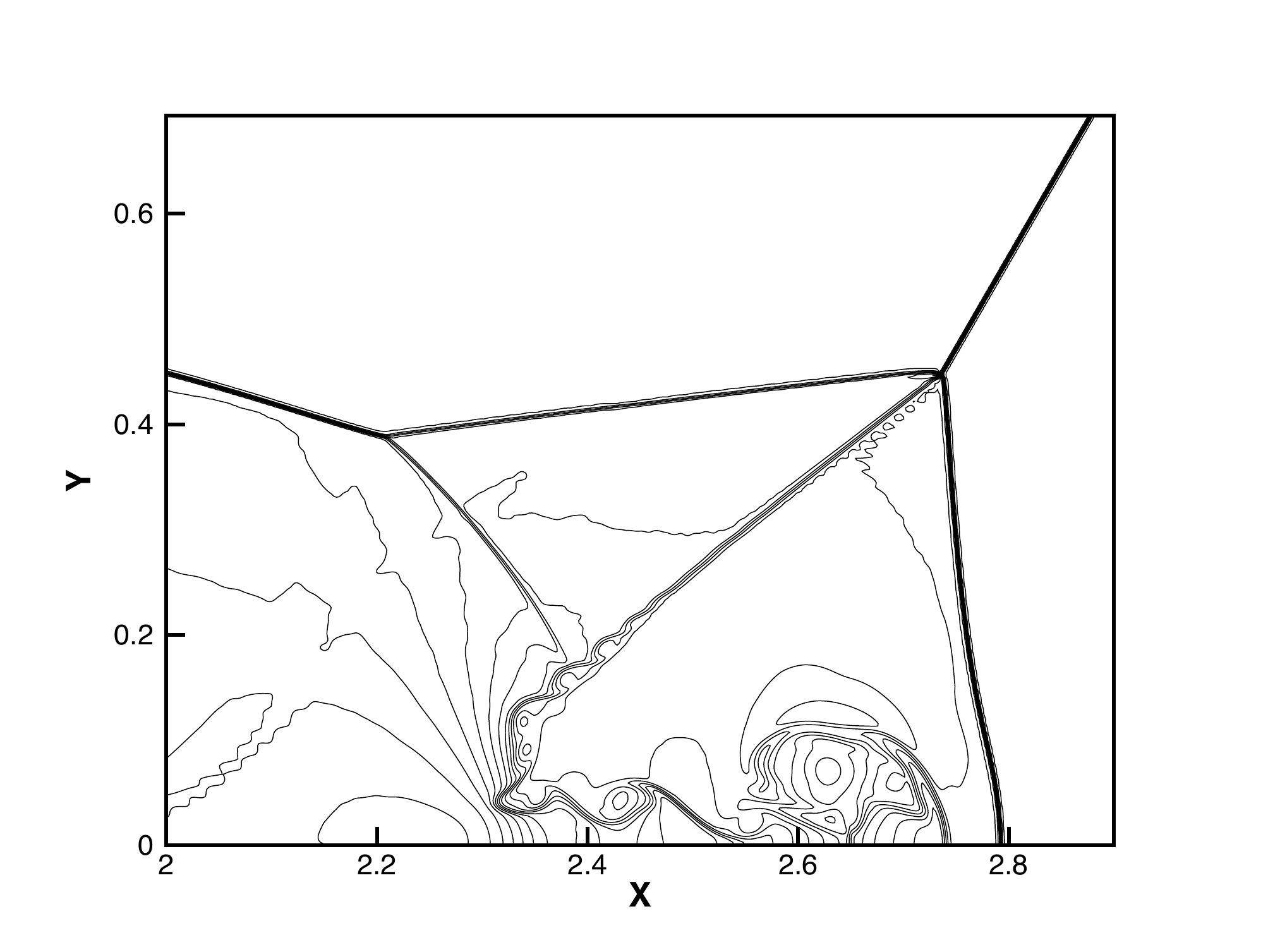,width=4 in} }
 \caption{ Double Mach reflection problem. T=0.2. 30 equally spaced density contours from 1.5 to
 22.7. From top to bottom: the results of the hybrid HWENO scheme; the cells where we modify the first order moments in the hybrid HWENO scheme at the final time; zoomed of the hybrid HWENO scheme. Uniform meshes with 1920 $\times$ 480 cells.}
\label{smhfig}
\end{figure}

\noindent {\bf Example 3.12.} Double Mach reflection problem \cite{Wooc}, which is modeled by the two-dimensional Euler equations (\ref{euler2}). The computational domain is
\( [0,4]\times[0,1]\). A reflection wall lies at the bottom,  starting from $x=\frac{1}{6}$, $y=0$, making a $60^{o}$ angle with the $x$-axis. For the bottom boundary, the exact post-shock  condition is imposed for the part from $x=0$ to $x=\frac{1}{6}$ and the rest is reflection boundary condition, while it is  the exact motion of the Mach 10 shock for the top boundary. $\gamma=1.4$ and the final computing time is set as $T=0.2$. In Figure \ref{smhfig}, we present the numerical results in region \( [0,3]\times[0,1]\), the cells where we modify the first order moments in the hybrid HWENO scheme at the final time and the blow-up region around the double Mach stems. Again, the hybrid HWENO scheme works well for this  test case, meanwhile, there are only 3.55\% cells  in which we need to modify their first order moments, which saves about 68.2\% computational time as most regions directly use linear approximation, and it shows the hybrid HWENO scheme has higher efficiency than HWENO scheme.
\smallskip

\begin{figure}[!ht]
\centerline{\psfig{file=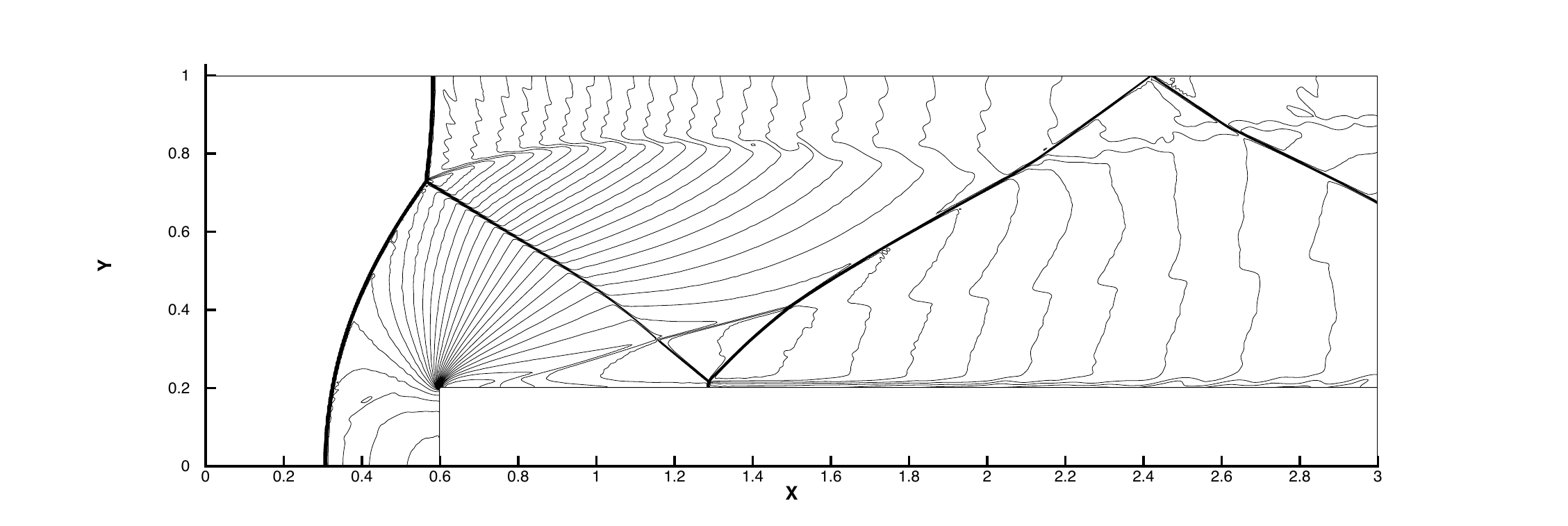,width=6 in}}
\centerline{\psfig{file=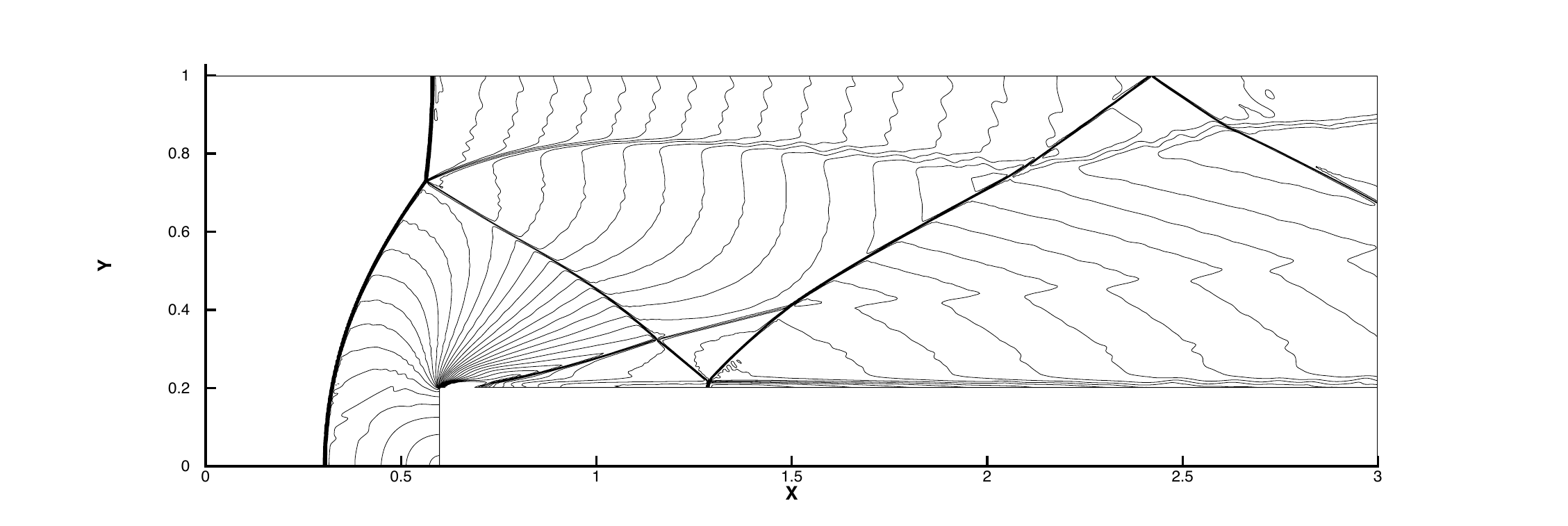,width=6 in}}
\centerline{\psfig{file=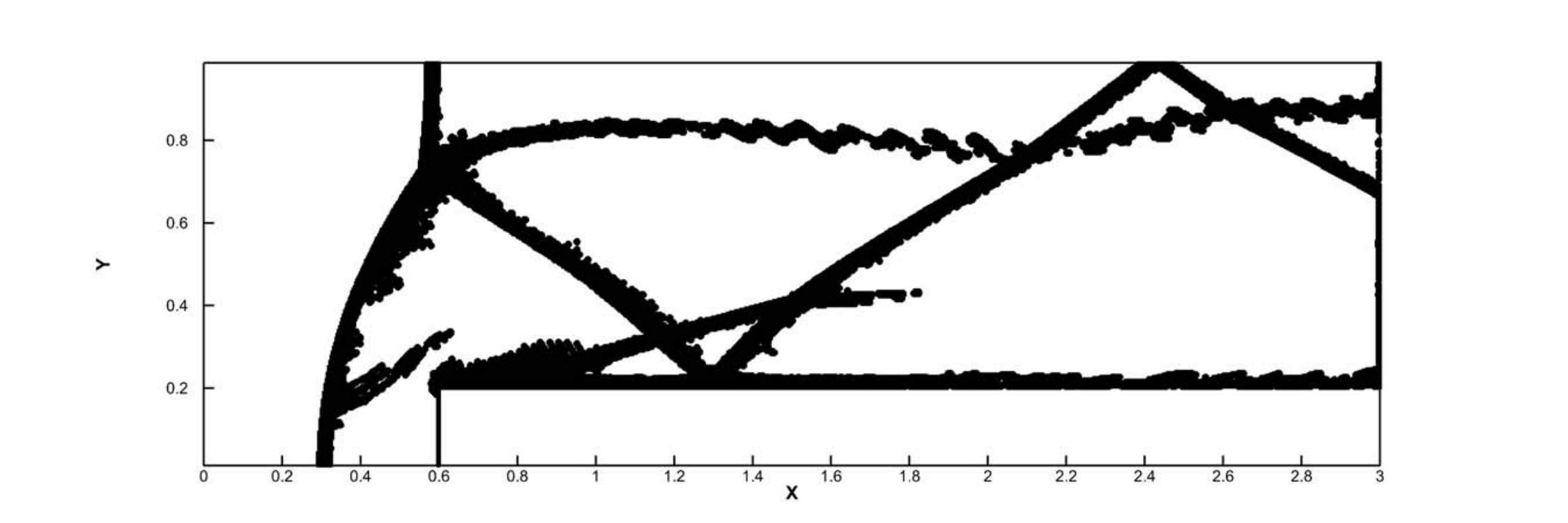,width=6 in}}
 \caption{ Forward step problem. T=4. From top to bottom: 30 equally spaced density contours from 0.32 to 6.15 of the hybrid HWENO scheme; 30 equally spaced Mach number contours from 0.05 to 3.05 of the hybrid HWENO scheme; the cells where we modify the first order moments in the hybrid HWENO scheme at the final time. Uniform meshes with 960 $\times$ 320 cells.}
\label{stepfig}
\end{figure}

\noindent {\bf Example 3.13.} The final example is about a Mach 3 wind tunnel with a step modeled by the two-dimensional Euler equations (\ref{euler2}), which is also originally from \cite{Wooc}. The setup of the problem is as follows. The wind tunnel is 1 length unit wide and 3 length units long. The step is 0.2 length units high and is located 0.6 length units from the left-hand end of the tunnel. The problem is initialized by a right-going Mach 3 flow. Reflective boundary conditions are applied along the wall of the tunnel and inflow/outflow  boundary conditions are applied at the entrance/exit. The corner of the step is a singular point and we treat it as the same way  in \cite{Wooc}, which is based on the assumption of a nearly steady flow in the region near the corner. The final time is $T=4$. In Figure \ref{stepfig}, we present the results for the hybrid HWENO scheme with $960\times320$ uniform cells and the cells where we modify the first order moments in the hybrid HWENO scheme at the last time step. We can notice that the hybrid  HWENO scheme gets good resolutions in the non-smooth  regions, moreover, we find that there are only 11.68\% cells where we need to modify their first order moments in our implementation, which shows most regions directly use high order linear approximation and we don't need to modify their first order moments, and we also know the hybrid HWENO scheme has higher efficiency than the HWENO scheme for saving near 64.2\% CPU time by calculating.
\smallskip

\section{Concluding remarks}
\label{sec4}
\setcounter{equation}{0}
\setcounter{figure}{0}
\setcounter{table}{0}

In this paper, a  hybrid  finite volume Hermite weighted essentially non-oscillatory (HWENO) scheme is designed for solving  hyperbolic conservation laws. Compared with other HWENO methods \cite{QSHw1,QSHw,ZQHW,C0,HXA,TQ,LLQ,ZA,TLQ,CZQ,LLQ2} for hyperbolic conservation laws, we bring the idea of limiter for discontinuous Galerkin (DG) method to modify the first order moments in the troubled cells, meanwhile, we find that there are only a small part of cells where we need to modify their first order moments, and we also reconstruct the interface point values for the solutions using HWENO approximation in these troubled cells. In other words, there are many cells in which we reconstruct the point values for the solutions directly employing linear approximation, which makes the hybrid HWENO scheme be higher efficiency. In general, the modification for the first order moments in the troubled cells is significant to avoid oscillations and keep  the scheme be robust for non-smooth numerical tests, and these smooth and non-smooth numerical results all illustrate  the good performances of the hybrid HWENO scheme.


\end{document}